%% file: LagrangeFEM_RD_2D.tex
\pdfoutput=1

\documentclass{siamart0216}

\usepackage{latexsym}
\usepackage{amssymb}
\usepackage{amsfonts}
\usepackage{amsmath}
\usepackage{epsfig}
\usepackage{graphicx}
\usepackage{layout}
\usepackage{indentfirst}
\usepackage{longtable}
\usepackage{multicol}
\usepackage{subfigure}
\usepackage{multirow}
\usepackage{wrapfig}
\usepackage{array}
\usepackage{tabularx}

\makeatletter
\def\blfootnote{\xdef\@thefnmark{}\@footnotetext}
\makeatother

\newcommand{\R}{\mathbb{R}}

\newcommand{\M}{\mathbf{M}}

\renewcommand{\P}{\mathbb{P}}

\newcommand{\cT}{{\cal T}}

\renewcommand{\d}{\hspace{1pt} {\rm d}}

\newcommand{\diff}[2]{\dfrac{d{#1}}{d{#2}}}

\newcommand{\dpar}[2]{\dfrac{\partial #1}{\partial #2}}

\newcommand{\bs}[1]{\boldsymbol{#1}}

\renewcommand{\leq}{\leqslant}
\renewcommand{\geq}{\geqslant}

\newcommand{\xnabla}{\nabla_{\bx}}
\newcommand{\be}{\mathbf e}
\newcommand{\bn}{\mathbf n}
\newcommand{\bx}{\mathbf{x}}

\newcommand{\bW}{\mathbf{W}}
\newcommand{\bX}{\mathbf{X}}
\newcommand{\bv}{\mathbf{v}}
\newcommand{\bvv}{\mathbf{u}}

\newcommand{\bt}{\mathbf{\theta}}
\newcommand{\bu}{\mathbf{u}}
\newcommand{\VV}{\mathcal{V}}
\newcommand{\EE}{\mathcal{E}}
\newcommand{\BPhi}{\mathbf{\Phi}}

\newcommand{\J}{\mathbb J}

\newcommand{\sigmav}{i_{\mathcal{V}}}
\newcommand{\sigmat}{i_{\mathcal{E}}}
\newcommand{\jsigmav}{j_{\mathcal{V}}}
\newcommand{\jsigmat}{j_{\mathcal{E}}}
\newcommand{\DV}{\mathcal{D}_{\mathcal{V}}}
\newcommand{\DE}{\mathcal{D}_{\mathcal{E}}}

\newtheorem{thm}{Theorem}[section]

\newtheorem{rem}[thm]{Remark}

\begin{document}

\newcommand{\TheTitle}{Multidimensional staggered grid residual distribution scheme for Lagrangian hydrodynamics} 
\newcommand{\TheAuthors}{R. Abgrall, K. Lipnikov, N. Morgan and S. Tokareva}

\headers{Multidimensional SGH RD scheme for Lagrangian hydrodynamics}{\TheAuthors}

\title{{\TheTitle}}

\author{
  R\'emi Abgrall\thanks{Institute of Mathematics, University of Zurich, Switzerland (\email{remi.abgrall@math.uzh.ch})}
  \and
  Konstantin Lipnikov\thanks{Los Alamos National Laboratory, USA (\email{lipnikov@lanl.gov})}
  \and
  Nathaniel Morgan\thanks{Los Alamos National Laboratory, USA (\email{nmorgan@lanl.gov})}
  \and
  Svetlana Tokareva\thanks{\emph{corresponding author}, Los Alamos National Laboratory, USA (\email{tokareva@lanl.gov})}
}

\maketitle

\begin{abstract}
	
We present the second-order multidimensional Staggered Grid Hydrodynamics Residual Distribution (SGH RD) scheme for Lagrangian hydrodynamics. The SGH RD scheme is based on the staggered finite element discretizations as in [Dobrev et al., SISC, 2012]. However, the advantage of the residual formulation over classical FEM approaches consists in the natural mass matrix diagonalization which allows one to avoid the solution of the linear system with the global sparse mass matrix while retaining the desired order of accuracy. This is achieved by using Bernstein polynomials as finite element shape functions and coupling the space discretization with the deferred correction type timestepping method. Moreover, it can be shown that for the Lagrangian formulation written in non-conservative form, our residual distribution scheme ensures the exact conservation of the mass, momentum and total energy. In this paper we also discuss construction of numerical viscosity approximations for the SGH RD scheme allowing to reduce the dissipation of the numerical solution. Thanks to the generic formulation of the staggered grid residual distribution scheme, it can be directly applied to both single- and multimaterial and multiphase models. Finally, we  demonstrate computational results obtained with the proposed residual distribution scheme for several challenging test problems.

\end{abstract}



\section{Introduction}
\label{sec:intro}

The present paper extends the results of \cite{SISC2017} to the multidimensional case. We are interested in the numerical solution of the Euler equations in Lagrangian form. It is well known that there are two formulations of the fluid mechanics equations, depending on whether the formulation is done in a fixed frame (Eulerian formulation) or a reference frame moving at the fluid speed (Lagrangian formulation). There is also an intermediate formulation, the ALE (Arbitrary Eulerian Lagrangian) formulation where the reference frame is moving at a speed that is generally different from the fluid velocity. Each of these formulations has advantages and drawbacks. The Eulerian one is conceptually the simplest because the reference frame is not moving; this implies that the computations are performed on a fixed grid. The two others are conceptually more complicated because of a moving reference frame; which means that the grid is moving and the mesh elements are changing shape and thus tangling of the elements is possible. 

However, moving reference frame is advantageous for computing shock waves, slip lines, contact discontinuities and material interfaces. Usually slip lines are difficult to compute because of excessive numerical dissipation, and hence dealing with a mesh that moves with the flow is a straightforward way to minimize this dissipation because the slip lines are steady in the Lagrangian frame. This nice property of a relatively simple and efficient way to deal with slip lines has motivated  many researchers, starting from the seminal work of von Neumann and Richtmyer \cite{vonNeuman}, to more recent works such as \cite{Shashkov,Maire,Barlow,scovazzi,Cheng-Shu-JCP2014,ChengShu2014,MARS2014,LANLDG1}.

Most of these works deal with schemes that are formally second order accurate. Up to our knowledge, there are much less works dealing with (formally) high order methods: either they are of discontinuous Galerkin type \cite{vilar1,vilar2,vilar3}, use a staggered finite element formulation \cite{rieben} or an ENO/WENO formalism \cite{Cheng-Shu-JCP2014}, see also the recent developments in \cite{Dumbser-2013-1,Dumbser-2014-1,Dumbser-2014-2,Dumbser-2015-2}.

In the discontinuous Galerkin (DG) formulation, all variables are associated to the elements, while in the staggered grid formulation, the approximations of the thermodynamic variables (such as pressure, specific internal energy or specific volume/density) are cell-centered, and thus possibly discontinuous across elements as in the DG method, while the velocity approximation is node-based, that is, it is described by a function that is polynomial in each element and globally continuous in the whole computational domain. In a way this is a natural extension of the Wilkins' scheme \cite{Wilkins} to higher order of accuracy.

This paper follows the finite element staggered grid approach of Dobrev et al. \cite{rieben}. This formulation involves two ingredients. First, the staggered discretization leads to a global mass matrix that is block diagonal on the thermodynamic parameters (as in DG method) and a sparse symmetric mass matrix for the velocity components (as in finite element method). Hence, the computations require the inversion\footnote{By saying "inversion of a matrix" we mean the solution of a linear system with the corresponding matrix.} of a block diagonal matrix, which is cheap, but also of a sparse symmetric positive definite matrix, which is more expensive both in terms of CPU time and memory requirements. In addition, every time when mesh refinement or remapping is needed (which is typical for Lagrangian methods), this global matrix needs to be recomputed. Second, an artificial viscosity technique is applied in order to make possible the computation of strong discontinuities. 

Our method relies on the Residual Distribution (RD) interpretation of the staggered grid scheme of \cite{rieben}, however the artificial viscosity term is introduced differently. See \cite{CAMWA18} and references therein for details about RD scheme for multidimensional Euler equations. The aims of this paper are the following: (i) extend the method of \cite{SISC2017} to two-dimensional staggered grid formulation avoiding the inversion of the large sparse global mass matrix while keeping all the accuracy properties and (ii) optimize the artificial viscosity term to provide low dissipation while retaining stability. We also present the way to ensure the conservation of the total energy, which is done similarly to \cite{SISC2017} and \cite{CAF18}.

The structure of this paper is the following. In Section~\ref{sec:equations}, we derive the formulation of the Euler equations in Lagrangian form and then in Section~\ref{sec:staggered_grid_scheme} recall the staggered grid formulation for multiple spacial dimensions. Next, in Section~\ref{sec:rd_scheme}, we recall the RD formulation for time-dependent problems. In Section~\ref{sec:timestepping}, we explain the diagonalization of the global sparse mass matrix without the loss of accuracy: this is obtained by modifying the timestepping method by applying ideas coming from \cite{mario,AbgrallBacigaluppiTokareva2016,SISC2017,CAMWA18,CAF18}. In Section~\ref{sec:rd_hydro}, we explain how to adapt the RD framework to the equations of Lagrangian hydrodynamics. In Section~\ref{sec:conservation}, we show how the conservation of the total energy is ensured. In Section~\ref{sec:viscosity}, we recall the construction of MARS (Multidirectional Approximate Riemann Solution) artificial viscosity terms from \cite{MARS2014,MARS2017} and incorporate them in the first-order residuals so that the numerical viscosity depends on the direction of the flow which reduces the overall numerical dissipation. We demonstrate the robustness of the proposed scheme by considering several challenging two-dimensional test problems in Section~\ref{sec:numerical_results}. 

\section{Governing equations}
\label{sec:equations}

We consider a fluid domain $\Omega_0\subset\R^d$, $d=1,2,3$ that is deforming in time through the movement of the fluid, the deformed domain is denoted by $\Omega_t$. In what follows, $\bX$ denotes any point of $\Omega_0$, while $\bx$ denotes any point of $\Omega_t$, the domain obtained from $\Omega_0$ under deformation. We assume the existence of a one-to-one mapping $\Phi$ from $\Omega_0$ to $\Omega_t$ such that $\bx=\Phi(\bX,t)\in \Omega_t$ for any $\bX\in \Omega_0$. We will call $\bX$ the Lagrangian coordinates and $\bx$ the Eulerian ones.
The Lagrangian description corresponds to the one of an observer moving with the fluid. In particular, its velocity, which coincides with the fluid velocity, is given by:
\begin{equation}
\label{Eq:velocity}
\bvv(\bx,t)=\dfrac{d\bx}{dt}=\dpar{\Phi}{t}(\bX,t).
\end{equation}
We also introduce the deformation tensor $\J$ (Jacobian matrix),
\begin{equation}
\label{Eq:deformation}
\J(\bx, t)=\nabla_{\bX}{\Phi} (\bX, t) \text{ where } \bx=\Phi(\bX,t).
\end{equation}
Hereafter, the notation $\nabla_{\bX}$ corresponds to the differentiation with respect to Lagrangian coordinates, while $\nabla_{\bx}$ corresponds to the Eulerian ones.

It is well known that the equations describing the evolution of fluid particles are consequences of the conservation of mass, momentum and energy, 
as well as a technical relation, the Reynolds transport theorem. It states that for any scalar quantity $\alpha(\bx,t)$, we have:
\begin{equation}
\label{Eq:Transport:Reynolds}
\dfrac{d}{dt}\int_{\omega_t} \alpha(\bx,t)\,d\bx = \int_{\omega_t} \dpar{\alpha}{t}(\bx,t)\,d\bx +\int_{\partial\omega_t} \alpha(\bx,t) \bvv\cdot \bn\,d\sigma = \int_{\omega_t} \bigg( \dfrac{d\alpha}{dt} + \alpha\nabla_{\bx}\cdot\bvv \bigg)\,d\bx
\end{equation}
In this relation, the set $\omega_t$ is the image of any set $\omega_0\subset \Omega_0$ by $\Phi$, i.e. $\omega_t=\Phi(\omega_0,t)$, $d\sigma$ is the measure on the boundary of $\partial \omega_t$  and $\bn$ is the outward unit normal. The gradient operator is taken with respect to the Eulerian coordinates. 

The conservation of mass reads: for any $\omega_0\subset\Omega_0$, 
\[ \dfrac{d}{dt}\int_{\omega_t} \rho\,d\bx=0, \quad \omega_t=\Phi(\omega_0,t), \]
so that  we get, defining $J(\bx,t)=\det \J(\bx,t)$,
\begin{equation}
\label{Eq:Lagrange-mass-tau}
J(\bx,t) \rho(\bx,t)=\rho(\bX,0) := \rho_0(\bX).
\end{equation}


Newton's law states that the acceleration is equal to the sum of external forces, so that
\[
\dfrac{d}{dt}\int_{\omega_t} \rho \bvv\,d\bx=-\int_{\partial \omega_t} \bs{\tau} \cdot \bn\,d\sigma,
\]
where $\bs{\tau}$ is the stress tensor\footnote{Here, if $\mathbf{X}$ is a tensor and $\mathbf{y}$ is a vector, $\mathbf{X}\cdot \mathbf{y}$ is the usual matrix-vector multiplication.}. Here, we have  $\bs{\tau}=-p\; \mathbf{Id}_d$
where the pressure $p(\bx,t)$ is a thermodynamic characteristic of a fluid and in the simplest case a function of two independent thermodynamic parameters, for example the specific energy $\varepsilon$ and the density,
\begin{equation}
\label{Eq:EOS}p=p(\rho, \varepsilon).
\end{equation}

The total energy of a fluid particle is $\rho e=\rho \varepsilon+\frac12\rho \bvv^2$. Using the first principle of thermodynamics, the variation of energy is the sum of variations of heat and the work of the external forces.  Assuming an isolated system, we get
\[\dfrac{d}{dt}\int_{\omega_t}\rho(\varepsilon+\frac12 \bvv^2)\,d\bx=-\int_{\partial \omega_t} \big (\bs{\tau}\cdot  \bvv\big )\cdot \bn\,d\sigma,\]
that is, for fluids, 
\[\dfrac{d}{dt}\int_{\omega_t}\rho(\varepsilon+\frac12 \bvv^2)\,d\bx=-\int_{\partial \omega_t} p \bvv\cdot \bn\,d\sigma.\]

These integral relations lead to the following formulation of conservation laws in Lagrangian reference frame:
\begin{equation}
\label{Eq:Lagrange}
\begin{split}
\bvv(\bx,t)&=\dfrac{d\bx}{dt}, \quad \bx =\Phi(\bX,t)\\
J(\bx,t) \rho(\bx,t)&=\rho(\bX,0) := \rho_0(\bX),\\
\rho\dfrac{d\bvv}{dt}&+\nabla_{\bx} p=0\\
\rho \diff{\varepsilon}{t} &+ p\nabla_\bx \cdot \bvv=0.
\end{split}
\end{equation}
where $p=p(\rho, \varepsilon)$.

\section{Staggered grid formulation}
\label{sec:staggered_grid_scheme}

Here we briefly recall the main ideas of the staggered grid method used in \cite{rieben}.
A semi-discrete approximation of \eqref{Eq:Lagrange} is introduced such that the velocity field $\bvv$ and coordinate $\bx$ belong to a kinematic space $\mathcal{V}\subset \big( H^1(\Omega_0) \big)^d$, where $d$ is the spacial dimension; $\mathcal{V}$ has a basis denoted by $\{w_{\sigmav}\}_{\sigmav\in \mathcal{D}_{\mathcal{V}}}$, the set $\mathcal{D}_{\mathcal{V}}$ is the set of kinematic degrees of freedom (DOFs) with the total number of DOFs given by $\#\mathcal{D}_{\mathcal{V}}=N_{\mathcal{V}}$. 
The thermodynamic quantities such as the internal energy $\varepsilon$ and pressure $p$ are discretized in a thermodynamic space $\mathcal{E}\subset L^2(\Omega_0)$. As before, this space is finite dimensional, and its basis is 
$\{\phi_{\sigmat}\}_{\sigmat\in \mathcal{D}_{\mathcal{E}}}$. The set $\mathcal{D}_{\mathcal{E}}$ is the set of thermodynamical degrees of freedom with the total number of DOFs $\#\mathcal{D}_{\mathcal{E}}=N_{\mathcal{E}}$. In the following, the subscript ${\mathcal{V}}$ (resp. ${\mathcal{E}}$) refers to kinematic (resp. thermodynamic) degrees of freedom.

The fluid particle position $\bx$ is approximated by:
\begin{subequations}\label{Eq:LagrangeDisc}
	\begin{equation}
	\label{Eq:R:position}
	\bx=\Phi(\bX,t)=\sum\limits_{{\sigmav}\in \DV} \bx_{{\sigmav}}(t)w_{{\sigmav}}(\bX).
	\end{equation}
	
	The domain at time $t$ is then defined by
	\[ \Omega_t=\{ \bx\in \R^d \text{ such that there exists }\,\bX\in \Omega_0 : \bx=\Phi(\bX,t)\} \]
	where $\Phi$ is given by \eqref{Eq:R:position}.
	
	The velocity field is approximated by:
	\begin{equation}
	\label{Eq:R:velocity}
	\bvv(\bx,t)=\sum\limits_{{\sigmav}\in \DV} \bvv_{{\sigmav}}(t) w_{{\sigmav}}(\bX),
	\end{equation}
	and the specific internal energy is given by:
	\begin{equation}
	\label{Eq:R:energy}
	\varepsilon(\bx,t)=\sum_{{\sigmat}\in \DE} \varepsilon_{{\sigmat}}(t) \phi_{{\sigmat}}(\bX).
	\end{equation}
	
	Considering the weak formulation of \eqref{Eq:Lagrange}, we get:
	\begin{enumerate}
		\item For the velocity equation, for any ${\sigmav} \in \mathcal{D}_{\mathcal{V}}$,  denoting  by $\bn$ the outward pointing unit vector of $\partial \Omega_t$,
		\begin{equation}
		\label{Eq:R:velocityDisc}
		\int_{\Omega_t} \rho \dfrac{d\bvv}{dt} w_{\sigmav}\,d\bx =-\int_{\Omega_t} \bs{\tau}:\xnabla  w_{\sigmav}\,d\bx+\int_{\partial \Omega_t}
		\big (\bs{\tau}\cdot\bn\big ) w_{\sigmav} \,d\sigma .
		\end{equation}

		Using \eqref{Eq:R:velocity}, we get\footnote{Here, if $\mathbf{X}$ and $\mathbf{Y}$ are tensors,  $\mathbf{X}:\mathbf{Y}$ is the contraction $\mathbf{X}:\mathbf{Y}=\text{trace}(\mathbf{X}^T\mathbf{Y})$.}
		\begin{equation*}
		\sum\limits_{j_{\VV}} \bigg( \int_{\Omega_t} \rho w_{j_\VV} w_{i_\VV}\,d\bx \bigg) \dfrac{d\bvv_{j_\VV}}{dt} = - \int_{\Omega_t} \bs{\tau}:\xnabla  w_{i_\VV}\,d\bx + \int_{\partial \Omega_t} \big (\bs{\tau}\cdot\bn\big ) w_{\sigmav}\,d\sigma.
		\end{equation*}
		
		Introducing the vector $\hat{\bvv}$ with components $\bvv_{i_\VV}$ and $\mathbf{F}$ the force vector given by the right-hand side of the above equation, we get the formulation
		\[ \M_{\mathcal{V}}\dfrac{d\hat{\bvv}}{dt} = \mathbf{F}. \]
		
		The kinematic mass matrix $\M_{\mathcal{V}}=(M_{i_\VV j_\VV}^{\mathcal{V}})$ has components 
		\[ M_{i_\VV j_\VV}^{\mathcal{V}}=\int_{\Omega_t} \rho w_{j_\VV} w_{i_\VV}\,d\bx. \]
		Thanks to the Reynolds transport theorem \eqref{Eq:Transport:Reynolds} and mass conservation, $\M_{\mathcal{V}}$ does not depend on time, see \cite{rieben} for details. Note that $\M_{\mathcal{V}}$ is a global, typically irreducible, sparse symmetric matrix of dimension $N_{\mathcal{V}}\times N_{\mathcal{V}}$ because the shape functions of $\DV$ are continuous.
		
		\item For the internal energy, we get a similar form,
		\begin{equation}
		\label{Eq:R:energyDisc}
		\int_{\Omega_t} \rho \dfrac{d\varepsilon}{dt} \phi_{i_\EE}\,d\bx = \int_{\Omega_t} \phi_{i_\EE} \bs{\tau}:\xnabla\bvv\,d\bx,
		\end{equation}
		which leads to 
		\[ \M_{\mathcal{E}} \dfrac{d\mathbf{\hat{\bs{\varepsilon}}}}{dt} = \mathbf{W}, \]
		where $\hat{\bs{\varepsilon}}$ is the vector with components $\varepsilon_{\sigmat}$, the thermodynamic mass matrix $\M_{\mathcal{E}}=(M_{i_\EE j_\EE}^{\mathcal{E}})$ with entries $M_{i_\EE j_\EE}^{\mathcal{E}}=\int_{\Omega_t} \phi_{i_\EE}\phi_{j_\EE}\,d\bx$ is again independent of time and $\bW$ is the right-hand side of \eqref{Eq:R:energyDisc}. Note that the thermodynamic mass matrix can be made block-diagonal by considering the shape functions with local support in $K \in \Omega_0$.
		
		\item The mass satisfies: 
		\begin{equation}
		\label{Eq:R:massDisc}
		\det{\J(\bx,t)} \rho(\bx, t)=\rho_0(\bX)
		\end{equation} where $\rho_0\in \mathcal{E}$ and the deformation tensor $\J$ is evaluated according to \eqref{Eq:R:position}.
		
		\item The positions $\bx_{{\sigmav}}$ satisfy:
		\begin{equation}
		\label{Eq:R:positions}
		\dfrac{d\bx_{{\sigmav}}}{dt}=\bvv_{{\sigmav}}(\bx_{{\sigmav}},t)
		\end{equation}
	\end{enumerate}
\end{subequations}

It  remains to define the discrete spaces $\mathcal{V}$ and $\mathcal{E}$. To do this,
we consider a \emph{conformal} triangulation of the initial computational domain $\Omega_0\subset\R^d$, $d=1,2,3$, which we shall denote by $\mathcal{T}_h$. We denote by $K$ any element of  $\mathcal{T}_h$ and assume for simplicity that 
$\cup_K K=\Omega_0$. The set of boundary faces is denoted by $\mathcal{B}$ and a generic boundary face is denoted by $f$, thus $\cup_{f\in \mathcal{B}}f=\partial\Omega_0$. As usual, denoting by $\P^r(K)$ the set of polynomials of degree at most $r$ defined on $K$, we consider two functional spaces (with integer $ r \geq 1$):
\[ \mathcal{V}=\{\bv \in L^2(\Omega_0)^d, \forall K, \bv_{|K}\in \P^r(K)^d\}\cap C^0(\Omega_0) \]
and
\[ \mathcal{E}=\{\bt\in L^2(\Omega_0), \forall K, \bt_{|K}\in \P^{r-1}(K)\}. \]

The matrix $\M_{\mathcal{E}}$ is symmetric positive definite block-diagonal while 
$\M_{\mathcal{V}}$ is only a \emph{sparse} symmetric positive definite matrix. 

The fundamental assumption made here is that the mapping $\Phi$ is bijective. In numerical situations, this can be hard to achieve for long-time simulations, and thus mesh remapping and re-computation of the matrices $\M_\mathcal{E}$ and $\M_\mathcal{V}$ must be done from time to time; this issue is however outside of the scope of this paper, see \cite{JCP50:review} for detailed discussion. 

The scheme defined by \eqref{Eq:LagrangeDisc} is linearly stable because of the choices of the test and trial functions, but only linearly stable. Since we are looking for possibly discontinuous solutions, one possible approach to ensure stability is to add mechanism of artificial viscosity \cite{vonNeuman,rieben} . The idea  amounts to  modifying the stress tensor $\bs{\tau}=-p\mathbf{Id}_d$ by $\bs{\tau}=-p\mathbf{Id}_d + \bs{\tau}_a(\bx, t)$, where the term $\bs{\tau}_a(\bx, t)$ specifies the artificial viscosity. We refer to \cite{rieben} for details on the construction of $\bs{\tau}_a(\bx, t)$.

It is possible to rewrite the system \eqref{Eq:Lagrange}, and in particular the relations \eqref{Eq:R:velocityDisc} and \eqref{Eq:R:energyDisc} in a slightly different way.
Let $K$ be any element of the triangulation $\mathcal{T}_h$, and for the kinematic degrees of freedom $\sigmav \in \DV$ and the thermodynamic degrees of freedom $\sigmat \in \DE$ consider the quantities 
\begin{align*}
\BPhi_{\VV,i_\VV}^K &= \int_K\bs{\tau}: \xnabla w_{i_\VV}\,d\bx -\int_{\partial K}\hat{\bs{\tau}}_\bn w_{i_\VV} d\sigma, \\
\BPhi_{\EE,i_\EE}^K &= -\int_K \phi_{i_\EE}^K \bs{\tau}: \xnabla\bvv\,d\bx,
\end{align*}
where $\hat{\bs{\tau}}_\bn$ is any numerical flux consistent with $\bs{\tau}\cdot \bn$, see e.~g. \cite{ToroBook}.

Using the compactness of the support of the basis functions $w_{i_\VV}$ and $\phi_{i_\EE}$, we can rewrite the relations \eqref{Eq:R:velocityDisc} and \eqref{Eq:R:energyDisc} as follows\footnote{Hereafter, we use the notation $\sum\limits_{K \ni i}$ to indicate that the summation is done over all elements $K$ containing a degree of freedom $i$}:

\begin{subequations}
	\label{RD:disc}
	\begin{equation}
	\label{RD:velocityDisc}
	\int_{\Omega_t} \rho \dfrac{d\bvv}{dt} w_{i_\VV}\,d\bx + \sum_{K \ni i_\VV} \BPhi_{\VV, i_\VV}^K = 0
	\end{equation}
	and 
	\begin{equation}
	\label{RD:energyDisc}
	\int_{\Omega_t} \rho \dfrac{d\varepsilon}{dt} \phi_{i_\EE}\,d\bx +\sum_{K \ni i_\EE} \BPhi_{\EE,i_\EE}^K=0,
	\end{equation}
	and we notice that on each element $K$, we have:
	\begin{equation}
	\begin{split}
	\label{Eq:conservation:R:E}
	\sum_{i_\EE\in K}&\BPhi_{\EE,i_\EE}^K + \sum_{i_\VV\in K} \bvv_{i_\VV} \cdot \BPhi_{\VV, i_\VV}^K
	=
	-\sum_{i_\EE\in K}\int_K \phi_{i_\EE}^K \bs{\tau}: \xnabla\bvv\,d\bx \\&+\sum_{i_\VV\in K}\bigg( \bvv_{i_\VV} \cdot \int_K \tau:\xnabla w_{i_\VV}\,d\bx - \bvv_{i_\VV}\cdot\int_{\partial K} \hat{\bs{\tau}}_\bn w_{i_\VV}\,d\sigma \bigg)\\
	&=-\int_K\bs{\tau}: \xnabla \bvv\,d\bx + \int_K \bs{\tau}:\xnabla \bvv\,d\bx - \int_{\partial K}\hat{\bs{\tau}}_\bn\cdot \bvv\,d\sigma =-\int_{\partial K}\hat{\bs{\tau}}_\bn\cdot\bvv\,d\sigma.
	\end{split}\end{equation}
\end{subequations}

There is no ambiguity in the definition of the last integral in \eqref{Eq:conservation:R:E} because $\bvv$ is continuous across $\partial K$ and the numerical flux 
$\hat{\bs{\tau}}_\bn$ is well defined.



\section{Residual distribution scheme}
\label{sec:rd_scheme}

In this section, we briefly recall the concept of residual distribution schemes for the following problem in $\Omega\subset \R^d$:
\[
\dpar{u}{t}+\xnabla \cdot \mathbf{f}(u)=0
\]
with the initial condition $u(\bx,0)=u_0(\bx)$. For simplicity we assume that $u$ is a real-valued function. Again, we consider a triangulation $\cT_h$ of $\Omega$. We want to approximate 
$u$ in 
\[ V_h=\{ u \in L^2(\Omega), \text{ for any } K\in \cT_h, u_{|K}\in \P^r\}\cap C^0(\Omega).\]
The set $\{\varphi_i\}$ is a basis of $V_h$, and $u_i$ are such that 
\begin{equation}
\label{approx_uh}
u=\sum_i u_i \varphi_i.
\end{equation}
As usual, $h$ represents the maximal diameter of the element of $\cT_h$.
We use  the same notations as before, and here the index $i$ denotes a generic degree of freedom.

\subsection{Residual distribution framework for steady problems}

We start by the steady problem, 
\[\xnabla \cdot \mathbf{f}(u)=0\]
and omit, for the sake of simplicity, the boundary conditions, see \cite{abgrall:cons} for details. We consider schemes of the form: for any degree of freedom $i$, 
\begin{equation}
\label{RD}
\sum_{K \ni i} \Phi_i^K(u)=0.
\end{equation}
The residuals must satisfy the conservation relation: for any $K$,
\begin{equation}
\label{RD:cons}
\sum_{i\in K} \Phi_i^K(u)=\int_{\partial K} \mathbf{f}^h\cdot \bn\,d\sigma := \Phi^K(u).
\end{equation}
Here, $\mathbf{f}^h\cdot \bn$ is an $(r+1)$-th order approximation of $\mathbf{f}(u)\cdot \bn$. Given a sequence of meshes that are shape regular with  $h\rightarrow 0$, one can construct a sequence of solution.
In \cite{AbgrallRoe2001}, it is shown that ,  if (i) this sequence of solutions stays bounded in $L^\infty$, (ii) a sub-sequence of it  converges in $L^2(\Omega)$ towards a limit $u$ and (iii) the residuals are continuous with respect to $u$, then the conservation condition guaranties that $u$ is a weak solution of the problem.

A typical example of such residual is the Rusanov residual,
\begin{equation*}
\Phi_i^{K,\text{Rus}}(u)=-\int_K \xnabla \varphi_i \cdot \mathbf{f}^h\,d\bx+\int_{\partial K} \mathbf{f}^h\cdot \bn\,\varphi_i\,d\sigma +
\alpha_K (u_i-\bar{u}_K),
\end{equation*}
where 
\[ \bar{u}_K = \dfrac{1}{N_K}\sum_{j\in K} u_j \]
with $N_K$ being the number of degrees of freedom inside an element $K$ and 
	\[\alpha_K \geq |K|\; \max\limits_{i\in K}\max\limits_{\bx \in K}\; \rho\big(\nabla_{\mathbf{u}}\mathbf{f}(u)\cdot \nabla \varphi_i).\]
	Here, $\rho(A)$ is the spectral radius of the matrix $A$.

This residual can be rewritten as 
\[\Phi_i^{K,\text{Rus}}(u)=\sum_{j\in K} c_{ij}^K (u_i-u_j)\]
with $c_{ji}^K\geq 0$.
It is easy to see that using the Rusanov residual leads to very dissipative solutions, but the scheme is easily shown to be monotonicity preserving in the scalar case, see for example \cite{AbgrallRoe2001}. There is a systematic way of improving the accuracy. One can show \cite{AbgrallRoe2001} that if the residuals satisfy, for any degree of freedom $i$, 
\[\Phi_i^K(u_{ex}^h)=O(h^{k+d}),\]
where $u_{ex}$ is the exact solution of the steady problem,  $u_{ex}^h$ is an interpolation of order $k+1$ and $d$ is the dimension of the problem, then the scheme is formally of order $k+1$. It is shown in \cite{AbgrallRoe2001} how to achieve a high order of accuracy while keeping the monotonicity preserving property. A systematic way of achieving this is to set:
\begin{equation}\label{Phi:psi}
\Phi_i^K(u)=\beta_i^K(u) \Phi^K(u),
\end{equation} where the distribution coefficients $\beta_i^K(u)$ are given by
\begin{equation}
\label{Psi2}
\beta_i^K(u)=\dfrac{\max\big( \frac{\Phi_i^{K,\text{Rus}}}{\Phi^K},0 \big)}{\sum\limits_{j\in K} \max\big( \frac{\Phi_j^{K,\text{Rus}}}{\Phi^K},0 \big)}
\end{equation}
and $\Phi^K$ is defined by \eqref{RD:cons}.
Some refinements exist in order to get an entropy inequality, see \cite{ABGRALL2018640} for example. Note that $\beta_i^K(u)$ is constant on $K$.

It is easy to see that one can rewrite \eqref{Phi:psi} in a Petrov-Galerkin fashion:
\begin{multline*}
\Phi_i^K(u) = \int_K \beta_i^K(u)\xnabla\cdot\mathbf{f}^h\,d\bx =  \int_K \varphi_i \xnabla\cdot\mathbf{f}^h\,d\bx+\int_K \big( \beta_i^K(u)-\varphi_i \big) \nabla_u \mathbf{f}^h\cdot \xnabla u\,d\bx\\
=-\int_K\xnabla \varphi_i \cdot \mathbf{f}^h\,d\bx +\int_{\partial K} \varphi_i\,\mathbf{f}^h \cdot \bn\,d\sigma +\int_K \big( \beta_i^K(u)-\varphi_i \big) \nabla_u  \mathbf{f}^h\cdot \xnabla u\,d\bx,
\end{multline*}
so that from \eqref{RD} we get
\[
0=-\int_\Omega \xnabla \varphi_i \cdot \mathbf{f}^h\,d\bx+\int_{\partial \Omega} \varphi_i \mathbf{f}^h(u)\cdot \bn\,d\sigma +\sum_{K \ni i} \int_K \big( \beta_i^K(u)-\varphi_i \big) \nabla_u \mathbf{f}^h\cdot \xnabla u\,d\bx.
\]

Inspired by this formulation, we would naturally discretize the unsteady problem as:
\begin{multline}
\label{RD:scheme_unsteady}
0=\int_\Omega\varphi_i \dpar{u}{t}\,d\bx-\int_\Omega \xnabla \varphi_i \cdot \mathbf{f}^h \,d\bx+\int_{\partial \Omega} \varphi_i\,\mathbf{f}^h(u)\cdot \bn\,d\sigma \\ +\sum_{K \ni i} \int_K \big( \beta_i^K(u)-\varphi_i \big) \bigg( \dpar{u}{t}+\nabla_u  \mathbf{f}^h\cdot \xnabla u \bigg)\,d\bx.
\end{multline}

The formulation \eqref{RD:scheme_unsteady} can be as well derived from \eqref{RD} by introducing the "space-time" residuals (the value of $\beta_i^K$ is not relevant at this stage)
\begin{equation}
\label{RD:Phi_unsteady}
\Phi_i^K(u)=\beta_i^K(u)\int_K\bigg( \dpar{u}{t}+\xnabla\cdot\mathbf{f}^h(u) \bigg)\,d\bx.
\end{equation}
The semi-discrete scheme \eqref{RD:scheme_unsteady} requires an appropriate ODE solver for time-stepping. 

A straightforward discretization of \eqref{RD:scheme_unsteady} would lead to a mass matrix $\M=(M_{ij})_{i,j}$ with entries
\[
M_{ij}=\int_\Omega\varphi_i \varphi_j\,d\bx + \sum_{K \ni  i}\int_K \big( \beta_i^K(u)-\varphi_i \big) \varphi_j\,d\bx.
\]
Unfortunately, this matrix has no special structure, might not be invertible (so the problem is not even well posed!), and in any case it is highly non linear since $\beta_i^K$ \emph{depends} on $u$.
A solution to circumvent the problem has been proposed in \cite{mario}. The main idea is to keep the spatial structure of the scheme and slightly modify the temporal one without violating the formal accuracy. A second order version of the method is designed in \cite{mario} and extension to high order is explained in \cite{AbgrallBacigaluppiTokareva2016}. For the purposes of this paper and for comparison with \cite{rieben} we only need the second order case.

\begin{figure}[H]
	\centering
	\subfigure[Step 1: Compute fluctuation]{\includegraphics[scale=0.3]{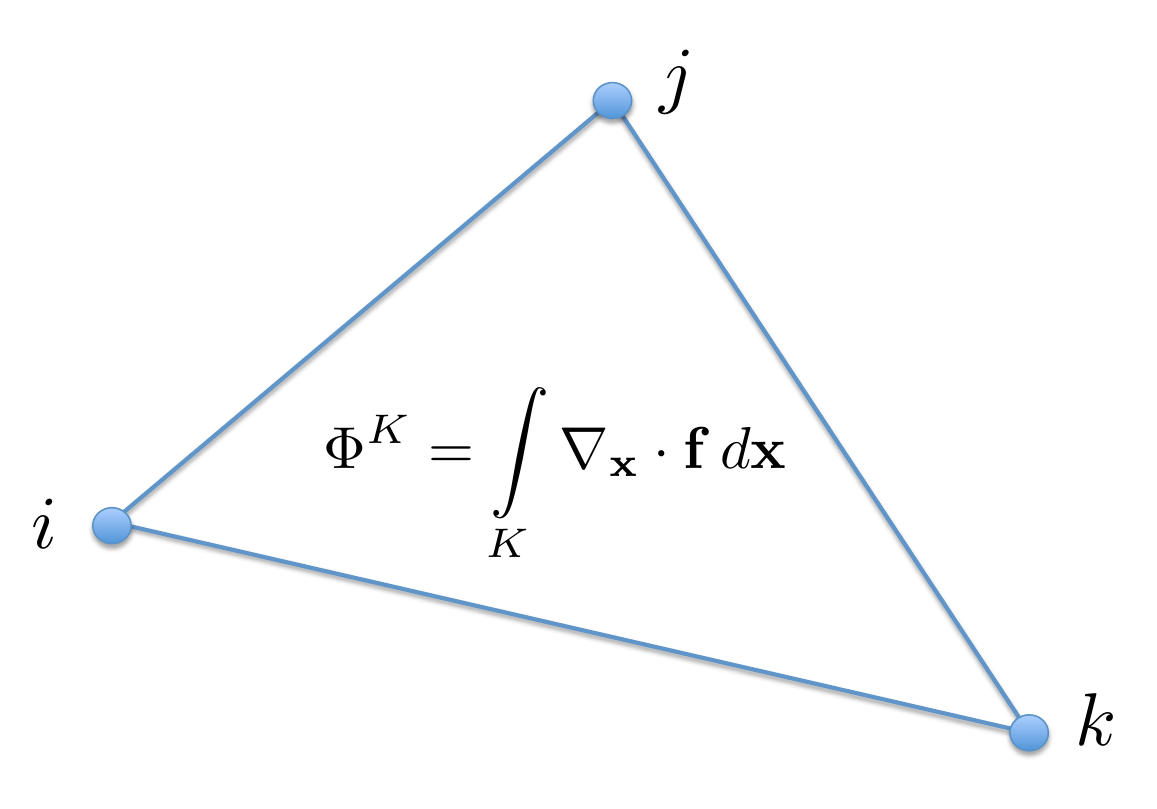}\label{Steps_RDa}}
	\hspace{0.18cm}\subfigure[Step 2: Split distribution]{\includegraphics[scale=0.3]{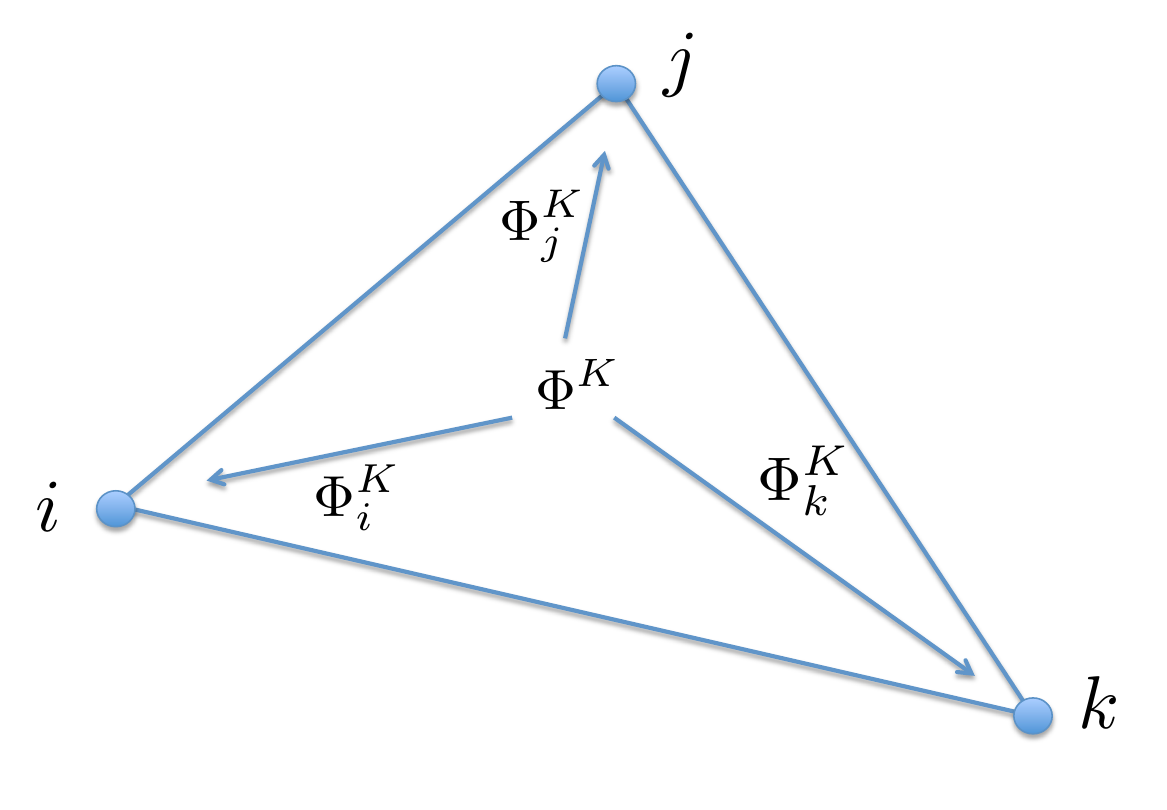}\label{Steps_RDb}}
	\hspace{0.33cm}\subfigure[Step 3: Gather residuals, evolve]{\includegraphics[scale=0.4]{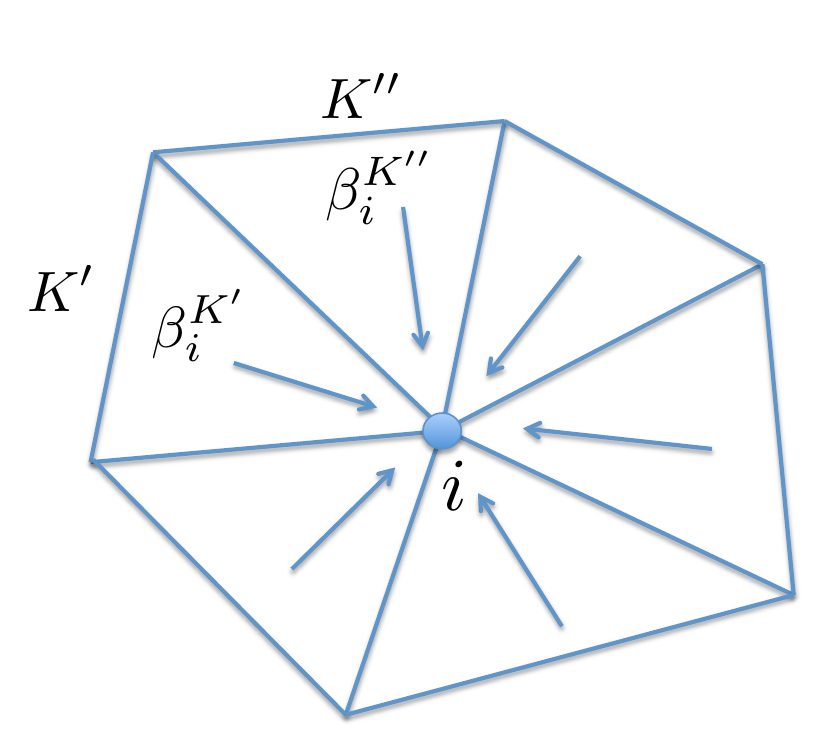}\label{Steps_RDc}}
	\caption{Illustration of the three steps of the residual distribution approach for linear triangular elements.}
	\label{Steps_RD}
\end{figure}

Hence, the main steps of the residual distribution approach could be summarized as follows (see also Fig.~\ref{Steps_RD} where the approach is illustrated for linear FEM on triangular elements):
\begin{enumerate}
	\item
	We define  for all $ K \in \Omega_h$ a fluctuation term (total residual), see Fig.~\ref{Steps_RDa} $$\Phi^K(u)=\int_K \nabla_x\cdot \mathbf{f}^h(u)\,d\mathbf{x}=\int_{\partial K} \mathbf{f}^h(u)\cdot\bn \; d\sigma$$ 
	
	\item
	We define a nodal residual $\Phi_i^K(u)$ as the contribution to the fluctuation term $\Phi^K$ from a degree of freedom (DOF) $i$ within the element $K$, so that the following conservation property holds (see Fig.~\ref{Steps_RDb}): for any element $K$ in $\Omega$,
	\begin{equation}
	\Phi^K(u)=\sum_{i \in K} \Phi_i^K(u), 
	\label{RD_distrib}
	\end{equation}

	The distribution strategy, i.e. how much of the fluctuation term has to be taken into account on each DOF $i \in K$, is defined by means of the distribution coefficients $\beta_i$:
	\begin{equation}
	\Phi_i^K=\beta_i^K \; \Phi^K,
	\label{def_phiK}
	\end{equation}
	where, due to \eqref{RD_distrib},
	\begin{equation*}
	\sum_{i \in K} \beta_i^K = 1.
	\end{equation*}
	
	\item The resulting scheme is obtained by collecting all the residual contributions $\Phi_i^K$ from elements $K$ surrounding a node $i \in \Omega$ (see Fig.~\ref{Steps_RDc}), that is
	\begin{equation}
	\sum_{K \ni i} \Phi_i^K(u) = 0, \quad \forall i \in \Omega,
	\label{RD_nodal}
	\end{equation}
	which allows to calculate the coefficients $u_i$ in the approximation \eqref{approx_uh}.
\end{enumerate}

\subsection{Second order timestepping method}
\label{sec:timestepping}

Here we describe the idea of the modified time stepping from \cite{mario}. We start with the description of our time-stepping algorithm based on a second order Runge-Kutta scheme for an ODE of the form 
\[y'+L(y)=0.\]

Given an approximate solution $y_n$ at time $t^n$, for the calculation of $y_{n+1}$ we proceed as follows:
\begin{enumerate}
	\item Set $y^{(0)}=y^n$;
	\item Compute $y^{(1)}$ defined by
	\[ \dfrac{y^{(1)}-y^{(0)}}{\Delta t}+L(y^{(0)})=0; \]
	\item Compute $y^{(2)}$ defined by
	\[ \dfrac{y^{(2)}-y^{(0)}}{\Delta t}+\dfrac{L(y^{(0)})+L(y^{(1)})}{2}=0; \]
	\item Set $y^{n+1}=y^{(2)}$.
\end{enumerate}

We see that the generic step in this scheme has the form
\[ \dfrac{\delta^ky}{\Delta t}+\mathcal{L}(y^{(0)},y^{(k)})=0\]
with 
\[ \mathcal{L}(a,b)=\dfrac{L(a)+L(b)}{2} \] 
and 
\[ \delta^k y=y^{(k+1)}-y^{(0)}, \quad k = 0,1.\]

A variant is to take $\mathcal{L}(a,b)=L\big ( \tfrac{a+b}{2}\big )$.

Coming back to the residuals \eqref{RD:Phi_unsteady}, we write for each element $K$ and $k=0,1$:
\begin{equation*}\begin{split}
\beta_i^K(u)\int_K&\bigg( \dfrac{\delta^k u}{\Delta t}+\mathcal{L}(u^{(0)},u^{(k)}) \bigg)\,d\bx  =
\int_K\varphi_i \bigg( \dfrac{\delta^k u}{\Delta t}
+\mathcal{L}(u^{(0)},u^{(k)}) \bigg)\,d\bx\\ &\qquad \qquad \qquad\qquad \qquad \qquad
+\int_K\big( \beta_i^K(u)-\varphi_i \big) \bigg( \dfrac{\delta^ku}{\Delta t}+\mathcal{L}(u^{(0)},u^{(k)}) \bigg)\,d\bx\\
&\approx \int_K\varphi_i \bigg( \dfrac{\delta^ku}{\Delta t}+\mathcal{L}(u^{(0)},u^{(k)}) \bigg)\,d\bx\\
&\qquad \qquad \qquad\qquad \qquad \qquad+ \int_K\big( \beta_i^K(u)-\varphi_i \big)\bigg( \dfrac{\widetilde{\delta^ku}}{\Delta t}+\mathcal{L}(u^{(0)},u^{(k)}) \bigg)\,d\bx
\end{split}\end{equation*}
with
\[ 
\widetilde{\delta^ku}=\left\{ \begin{array}{ll}
0 & \text{ if } k=0,\\
u^{(1)}-u^{(0)} & \text{ if } k=1.
\end{array}
\right.
\]
We see that
\begin{multline*}
\int_K\varphi_i \bigg( \dfrac{\delta^ku}{\Delta t}+\mathcal{L}(u^{(0)},u^{(k)}) \bigg)\,d\bx
+ \int_K\big( \beta_i^K(u)-\varphi_i \big) \bigg( \dfrac{\widetilde{\delta^ku}}{\Delta t}+\mathcal{L}(u^{(0)},u^{(k)}) \bigg)\,d\bx\\
=\int_K \varphi_i\bigg(\dfrac{\delta^ku}{\Delta t}-\dfrac{\widetilde{\delta^ku}}{\Delta t} \bigg)\,d\bx + \beta_i^K(u)\int_K \bigg( \dfrac{\widetilde{\delta^ku}}{\Delta t}+\mathcal{L}(u^{(0)},u^{(k)}) \bigg)\,d\bx.
\end{multline*}
This relation is further simplified if mass lumping can be applied: letting 
\begin{equation}
\label{Eq:lumping_K}
C_i^K = \sum_{j \in K}\int_K\varphi_i\varphi_j\,d\bx = \int_K\varphi_i\,d\bx
\end{equation}
and 
\begin{equation}
\label{Eq:lumping_Omega}
C_i = \int_\Omega\varphi_i\,d\bx = \sum_{K \ni  i} \int_K\varphi_i\,d\bx,
\end{equation}
for the degree of freedom $i$ and the element $K$ we look at the quantity
\[
C_i^K\bigg( \dfrac{\delta^ku_i}{\Delta t}-\dfrac{\widetilde{\delta^ku_i}}{\Delta t} \bigg)+
\beta_i^K(u)\int_K \bigg( \dfrac{\widetilde{\delta^ku}}{\Delta t}+\mathcal{L}(u^{(0)},u^{(k)}) \bigg)\,d\bx,
\]
i.e.
\[
C_i^K \dfrac{u_i^{(k+1)}-u_i^{(k)}}{\Delta t} +\beta_i^K(u)\int_K \bigg( \dfrac{\widetilde{\delta^ku}}{\Delta t}+\mathcal{L}(u^{(0)},u^{(k)}) \bigg)\,d\bx.
\]
Here $\beta_i^K(u)$ is evaluated using \eqref{Psi2} where $\Phi_i^{K,\text{Rus}}$ is replaced by the modified space-time Rusanov residuals
\[
\Phi_i^{K,\text{Rus}}=\int_K\varphi_i \bigg( \dfrac{\widetilde{\delta^ku}}{\Delta t}+\mathcal{L}(u^{(0)},u^{(k)}) \bigg)\,d\bx + \dfrac{1}{2}\Big( \alpha_K^{(0)}
\big(u^{(0)}_i-\bar{u}_K^{(0)}\big) + \alpha_K^{(k)}
\big(u^{(k)}_i-\bar{u}_K^{(k)}\big)\Big),
\]
where 
\[\bar{u}^{(l)}=\frac{1}{N_K}\sum\limits_{j\in K} u^{(l)}_j, \quad l = 0,k,\] 
with $N_K$ being the number of degrees of freedom in an element $K$
and $\alpha_K$ large enough and, finally, 
\[ \Phi^K(u) = \sum_{i\in K}\Phi_i^{K,\text{Rus}}. \]
Then the idea is to use \eqref{RD} at each step of the Runge-Kutta method with the residuals given by
\begin{multline}
\label{RD:RK2}
\Phi_i^K(u)=\int_K\varphi_i \bigg( \dfrac{\delta^ku}{\Delta t}+\mathcal{L}(u^{(0)},u^{(k)}) \bigg)\,d\bx \\ + \int_K\big( \beta_i^K(u)-\varphi_i \big) \bigg( \dfrac{\widetilde{\delta^ku}}{\Delta t}+\mathcal{L}(u^{(0)},u^{(k)}) \bigg)\,d\bx,
\end{multline}
so that the overall step writes: for $k=0,1$ and any $i$,
\begin{equation}
\label{RD:RK3}
C_i \dfrac{u_i^{(k+1)}-u_i^{(k)}}{\Delta t}+\sum_{K \ni i} \Phi_{i,ts}^{K,L}=0,
\end{equation}
where we have introduced the limited space-time residuals
\begin{equation}
\Phi_{i,ts}^{K,L} = \beta_i^K\int_K \bigg( \dfrac{\widetilde{\delta^ku}}{\Delta t}+\mathcal{L}(u^{(0)},u^{(k)}) \bigg)\,d\bx.
\end{equation}
One can easily see that each step of \eqref{RD:RK3} is purely explicit. 

One can show that this scheme is second order in time. The \emph{key} reason for this is that
we have
\[ \sum_{i\in K} \int_K \big( \varphi_i - \beta_i^K(u) \big)\,d\bx=0, \]
see \cite{mario,AbgrallLaratRicchiuto2011} for details. 

\begin{rem}
	We need that $C_i>0$ for any degree of freedom. This might not hold, for example, for quadratic Lagrange basis. For this reason, we will use Bernstein elements for the approximation of the solution.
\end{rem}

\section{Residual distribution scheme for Lagrangian hydrodynamics}\label{sec:rd_hydro}

In this section, we explain how to adapt the previously derived RD framework to the equations of Lagrangian hydrodynamics. 
We consider the same functional spaces as in section~\ref{sec:staggered_grid_scheme}, namely the kinematic space $\mathcal{V}$  and the thermodynamic space $\mathcal{E}$. 

In the case of a simplex $K\subset \R^d$, one can consider the barycentric coordinates associated to the vertices of $K$ and denoted by
$\{\Lambda_j\}_{j=1, d+1}$. By definition, the barycentric coordinates are positive on $K$ and we can consider the Bernstein polynomials of degree $r$: define $r=i_1+\ldots +i_{d+1}$, then
\begin{equation}
\label{Eq:Bezier}
B_{i_1 \ldots i_{d+1}}=\dfrac{r!}{i_1! \cdots i_{d+1}!} \Lambda_1^{i_1}\ldots \Lambda_{d+1}^{i_{d+1}}. 
\end{equation}
Clearly, $B_{i_1 \ldots i_{d+1}}\geq 0$ on $K$ and using the binomial identity
\[ \sum\limits_{i_1, \ldots, i_{d+1}, \sum_{1}^{d+1}i_j=r}B_{i_1 \ldots i_{d+1}}=\bigg( \sum_{i=1}^{d+1}\Lambda_i \bigg)^r=1. \]
In the case of quadrilateral/hexahedral elements, there exists a mapping that transforms this element into the unit square/cube. Then we can proceed by tensorization of segments $[0,1]$ seen as one-dimensional simplicies.

It is left to define the residuals for the equations of the Lagrangian hydrodynamics. Since the PDE on the velocity is written in conservation form, the derivations presented in the previous section can be directly applied, see also \cite{mario} for the multidimensional case. However, we need to introduce some modifications for the thermodynamics. 
To this end, we first focus on the spatial term, in the spirit of \cite{mario,AbgrallLaratRicchiuto2011}.  We construct a first order monotone scheme, and using the technique of \cite{AbgrallLaratRicchiuto2011}, we design a formally high order accurate scheme. Therefore, we introduce the total residuals
\begin{equation}
\label{Eq:residual-tot-velo}
\Phi^K = \int_{\partial K} \hat{p}_n\,d\sigma \quad \text{and} \quad \Psi^K = \int_K p \nabla_\bx\cdot \bvv\,d\bx,
\end{equation}
where $p\in \mathcal{E}$,  $u\in \mathcal{V}$ and $\hat{p}_n$ is a consistent numerical flux which depends on the left and right state at $\partial K$. Next, the Galerkin residuals are given by
\begin{equation}
\label{Eq:residual-velo}
\begin{split}
\Phi_{\sigmav}^K &= -\int_K p\nabla_{\bx} \phi_{\sigmav}\,d\bx+\int_{\partial K}\phi_{\sigmav}\hat{p}_n\,d\sigma, \\
\Psi_{\sigmat}^K &= \int_K \psi_{\sigmat} p\nabla_\bx\cdot \bvv\,d\bx.
\end{split}
\end{equation}

From \eqref{Eq:residual-velo}, we define the Rusanov residuals
\begin{equation}
\label{Eq:residual-alpha-velo}
\Phi_{\sigmav}^{K,\text{Rus}}(\bvv,\varepsilon) = \Phi_{\sigmav}^K (\bvv,\varepsilon)+ \alpha_{K}(\bvv_{\sigmav}-\bar{\bvv}), \quad \bar{\bvv} = \dfrac{1}{N_{\mathcal{V}}^K} \sum_{{\sigmav}\in K} \bvv_{\sigmav}
\end{equation}
and 
\begin{equation}
\label{Eq:residual-alpha-eps}
\Psi_{\sigmat}^{K,\text{Rus}}(\bvv,\varepsilon) = \Psi_{\sigmat}^K(\bvv,\varepsilon) + \alpha_K(\varepsilon_{\sigmat}-\bar{\varepsilon}), \quad \bar{\varepsilon} = \dfrac{1}{N_{\mathcal{E}}^K} \sum_{{\sigmat}\in K} \varepsilon_{\sigmat}
\end{equation}
where $\alpha_K$ is an upper bound of the Lagrangian speed of sound $\rho c$ on $K$ multiplied by the measure of $\partial K$, and $N_{\mathcal{V}}^K$ (resp. $N_{\mathcal{E}}^K$) is the number of degrees of freedom for the velocity (resp. energy) on $K$.

The temporal discretization is done using the technique developed in the previous section. We introduce the modified space-time Rusanov residuals, for $k=0,1$:
\begin{equation}
\label{Eq:residual-ts-velo}
\Phi_{{\sigmav},ts}^{K,\text{Rus}} = \int_K \varphi_{\sigmav} \rho \dfrac{\widetilde{\delta^k \bvv}}{\Delta t}\,d\bx + \frac{1}{2}\bigg( \Phi_{\sigmav}^{K,\text{Rus}}(\bvv^{(0)},\varepsilon^{(0)})+\Phi_{\sigmav}^{K,\text{Rus}}(\bvv^{(k)},\varepsilon^{(k)}) \bigg)
\end{equation}
and
\begin{equation}
\label{Eq:residual-ts-eps}
\Psi_{{\sigmat},ts}^{K,\text{Rus}} = \int_K \psi_{\sigmat} \rho \dfrac{\widetilde{\delta^k \varepsilon}}{\Delta t}\,d\bx + 
\frac{1}{2}\bigg( \Psi_{\sigmat}^{K,\text{Rus}}(\bvv^{(0)},\varepsilon^{(0)})+\Psi_{\sigmat}^{K,\text{Rus}}(\bvv^{(k)},\varepsilon^{(k)}) \bigg).
\end{equation}

Finally, the high-order limited residuals are computed similarly to \eqref{Phi:psi} as
\begin{equation}
\label{Eq:residual-lim}
\Phi_{{\sigmav},ts}^{K,L} = \beta^K_{\sigmav}\Phi^K_{ts}, \quad \Psi_{{\sigmat},ts}^{K,L} = \beta^K_{\sigmat}\Psi^K_{ts},
\end{equation}
where the space-time Rusanov residuals \eqref{Eq:residual-ts-velo} and \eqref{Eq:residual-ts-eps} are used in expressions analogous to \eqref{Psi2} to calculate $\beta^K_{\sigmav}$ and $\beta^K_{\sigmat}$, respectively:
\begin{equation}
\label{beta-v-ts}
\beta^K_{\sigmav} = \dfrac{\max\big( \frac{\Phi_{{\sigmav},ts}^{K,\text{Rus}}}{\Phi^K_{ts}},0 \big)}{\sum\limits_{{\jsigmav} \in K} \max\big( \frac{\Phi_{{\jsigmav},ts}^{K,\text{Rus}}}{\Phi^K_{ts}},0 \big)},
\end{equation}
\begin{equation}
\label{beta-eps-ts}
\beta^K_{\sigmat} = \dfrac{\max\big( \frac{\Psi_{{\sigmat},ts}^{K,\text{Rus}}}{\Psi^K_{ts}},0 \big)}{\sum\limits_{{\jsigmat} \in K} \max\big( \frac{\Psi_{{\jsigmat},ts}^{K,\text{Rus}}}{\Psi^K_{ts}},0 \big)},
\end{equation}
and 
\begin{align*}
& \Phi^K_{ts} = \sum_{\sigmav \in K}\Phi_{{\sigmav},ts}^{K,\text{Rus}} = \int_K \bigg( \rho \dfrac{\widetilde{\delta^k \bvv}}{\Delta t} + \dfrac12\Big(\xnabla p^{(0)} + \xnabla p^{(k)}\Big) \bigg)\,d\bx, \\
& \Psi^K_{ts} = \sum_{\sigmat \in K}\Psi_{{\sigmat},ts}^{K,\text{Rus}} = \int_K \bigg( \rho \dfrac{\widetilde{\delta^k \varepsilon}}{\Delta t} + \dfrac12\Big(p^{(0)}\xnabla\cdot\bvv^{(0)} + p^{(k)}\xnabla\cdot\bvv^{(k)}\Big) \bigg)\,d\bx.
\end{align*}

Next, we introduce
\[ C^{\mathcal{V},K}_{\sigmav}=\int_K\rho \varphi_{\sigmav}\,d\bx, \qquad C^{\mathcal{E},K}_{\sigmat }=\int_K\rho \psi_{\sigmat}\,d\bx. \]
After applying the mass lumping as in \eqref{Eq:lumping_K}, \eqref{Eq:lumping_Omega}, the mass matrices $\M_\mathcal{V}$ for the velocity and $\M_{\mathcal{E}}$ for the thermodynamics become diagonal with entries at the diagonals given by
\begin{align*}
C_{\sigmav}^\mathcal{V} & =\int_\Omega\rho \varphi_{\sigmav}\,d\bx = \sum_{K\ni \sigmav} \int_K\rho \varphi_{\sigmav}\,d\bx, \\
C^\mathcal{E}_{\sigmat} &= \int_\Omega\rho \psi_{\sigmat}\,d\bx = \sum_{K\ni \sigmat} \int_K\rho \psi_{\sigmat}\,d\bx.
\end{align*}

Both matrices are invertible because $\varphi_{\sigmav}>0$ and $\psi_{\sigmat} >0$ in the element since we are using Bernstein basis. Note that we could have omitted the mass lumping for the thermodynamic relation because the mass matrix is block diagonal.

By construction, the scheme is conservative for the velocity, however, nothing is guaranteed for the specific energy. In order to solve this issue, inspired by the calculations of section~\ref{sec:staggered_grid_scheme},  and given a set of velocity residuals $\{\Phi_{\sigmav}^K\}$ and internal energy residuals $\{\Psi_{\sigmat}^K\}$, we slightly modify the internal energy evaluation by defining, together with \eqref{Eq:residual-lim},
\begin{equation}
\label{Eq:residual-lim-eps-cons}
\Psi_{\sigmat,ts}^{K,c} = \Psi_{\sigmat,ts}^{K,L}+r_{\sigmat},
\end{equation}
where the correction term $r_{\sigmat}$ is chosen to ensure the discrete conservation properties and will be specified in the following section. 

With all the above definitions,  the resulting residual distribution scheme is written as follows: for $k=0,1$
\begin{subequations}
	\label{eq:RD_scheme}
	\begin{equation}
	\label{Eq:RD_scheme-velo}
	C_{\sigmav}^{\mathcal{V}}\dfrac{\bvv_{\sigmav}^{(k+1)} - \bvv_{\sigmav}^{(k)}}{\Delta t}  + \sum_{K \ni \sigmav}
	\Phi_{{\sigmav},ts}^{K,L} = 0, 
	\end{equation}
	\begin{equation}
	\label{Eq:RD_scheme-eps}
	C_{ \sigmat}^{\mathcal{E}}\dfrac{\varepsilon_{\sigmat}^{(k+1)}-\varepsilon_{\sigmat}^{(k)}}{\Delta t}  + \sum_{K \ni \sigmat}
	\Psi_{{\sigmat},ts}^{K,c} = 0,
	\end{equation}
	\begin{equation}
	\label{Eq:RD_scheme-coor}
	\dfrac{\bx_{\sigmav}^{(k+1)}-\bx_{\sigmav}^n}{\Delta t} = \dfrac12\Big( \bvv_{\sigmav}^n + \bvv_{\sigmav}^{(k)} \Big).
	\end{equation}
\end{subequations}
Note that the discretization \eqref{Eq:RD_scheme-coor} is nothing but a second-order SSP RK scheme.

\section{Discrete conservation}
\label{sec:conservation}

Here we derive the expression for the term $r_{\sigmat}$ to ensure the local conservation property of the residual distribution scheme \eqref{eq:RD_scheme} and then give some conditions on the discrete entropy production.


The continuous problem satisfies the following conservation property for the specific total energy $e = \frac12\bvv^2+\varepsilon$:
\begin{equation}
\label{Eq:conserv-cont}
\int\limits_K\rho\diff{e}{t}\,dx + \int\limits_{\partial K} p \bvv\cdot \bn\,d\sigma = 0.
\end{equation}
The numerical scheme has to satisfy a conservation property analogous to \eqref{Eq:conserv-cont} at the discrete level. To achieve this, the thermodynamic residual has been modified according to \eqref{Eq:residual-lim-eps-cons}. 

The term $r_{\sigmat}$ is chosen such that:
\begin{multline}
\label{Eq:conserv-discr}
\sum_{\sigmav \in K} \tilde{\bvv}_{\sigmav}\bigg(  C_{\sigmav}^{\mathcal{V},K} \bigg( \dfrac{\delta^k \bvv_{\sigmav}}{\Delta t} - \dfrac{\widetilde{\delta^k \bvv_{\sigmav}}}{\Delta t}  \bigg) + \Phi_{\sigmav}^{K,L} \bigg) \\+ \sum_{\sigmat \in K} \bigg( C_{\sigmat}^{\mathcal{E},K} \bigg( \dfrac{\delta^k \varepsilon_{\sigmat}}{\Delta t} - \dfrac{\widetilde{\delta^k \varepsilon_{\sigmat}}}{\Delta t} \bigg)+ \Psi_{\sigmat,ts}^{K,L} + r_{\sigmat} \bigg) = 0,
\end{multline}
where 
\[
\tilde{\bu}_{\sigmav} = \dfrac12\big(\bu_{\sigmav}^{(k)} + \bu_{\sigmav}^{(k+1)}\big).
\]
We justify this relation in Appendix \ref{appendix:conservation}.

Since we have only one constraint, we impose in addition that $r_{\sigmat} = r$ for any $\sigmat \in K$, so that from \eqref{Eq:conserv-discr} we can derive
\begin{multline}
\label{correction:conservation}
r_{\sigmat} = \frac{1}{N_{\mathcal{E}}^K}\bigg( \int_{\partial K} \hat{p}_n \bvv\,d\sigma - \sum_{\sigmav \in K} \tilde{\bvv}_{\sigmav}\Phi_{\sigmav}^{K,L}  +  \sum_{\sigmav \in K} C_{\sigmav}^{\mathcal{V},K} \dfrac{\widetilde{\delta^k \bvv_{\sigmav}}}{\Delta t}  \\ - \sum_{\sigmat \in K} \Psi_{\sigmat}^{K,L} + \sum_{\sigmat \in K} C_{ \sigmat}^{\mathcal{E},K} \dfrac{\widetilde{\delta^k \varepsilon_{\sigmat}}}{\Delta t} \bigg) ,
\end{multline}
where $\hat{p}_n$ is the approximation of the pressure flux $p\bn$ at the boundary of the element $K$.

So far, we have indicated a way to recover local conservation by adding a term to the internal energy equation. This term depends on the residuals that are themselves constructed from first order residuals and in turn depend on the pressure flux, so that the conservation property is valid for any pressure flux. It is possible to add further constraints for better conservation properties and in this section we show how to impose a local (semi-discrete) entropy inequality. We also state two results that are behind the construction.

The discrete entropy production is discussed in \cite{SISC2017}: note that it is generic and does not use the fact whether the problem is one-dimensional or multidimensional.

\subsection{Entropy balance}

Since at the continuous level
\[ T\dfrac{ds}{dt}=\dfrac{d\varepsilon}{dt}+p\dfrac{dv}{dt}, \]
where $v = 1/\rho$ is the specific volume, and knowing that 
\[ \rho \dfrac{dv}{dt}=\nabla_\bx\cdot \bvv, \]
we look at the entropy inequality
\begin{equation}
\label{entropy:1}
\int_K \rho T\dfrac{ds}{dt} = \int_K \rho \bigg( \dfrac{d\varepsilon}{dt}+p\dfrac{dv}{dt} \bigg)\,d\bx=\int_K \bigg( \rho \dfrac{d\varepsilon}{dt} + p\nabla_\bx \cdot\bvv \bigg)\,d\bx \geq 0
\end{equation}
and try to derive its discrete counterpart.

For the sake of simplicity we demonstrate the discrete entropy balance conditions on the \emph{first-order} version of the scheme \eqref{eq:RD_scheme}. Taking the sum over the degrees of freedom of an element $K$ in equation \eqref{Eq:RD_scheme-eps} and noting that in the first-order scheme $\widetilde{\delta^k\varepsilon}/\Delta t = 0$ and $\Psi_{\sigmat}^{K,L} = \Psi_{\sigmat}^{K,\text{Rus}}$, we get
\begin{multline}
\label{entropy:2}
\sum_{\sigmat \in K} C_{ \sigmat}^{\mathcal{E},K} \dfrac{\delta^k \varepsilon_{\sigmat}}{\Delta t}  + \sum_{\sigmat \in K} \Psi_{\sigmat}^{K,c} \\= \sum_{\sigmat \in K} \bigg(\int_K \rho\psi_{\sigmat}\,d\bx \bigg)\dfrac{\delta^k \varepsilon_{\sigmat}}{\Delta t} + \sum_{\sigmat \in K} \big( \Psi_{\sigmat}^{K,\text{Rus}} + r_{\sigmat} \big)  = \int_K\rho\dfrac{\delta^k \varepsilon}{\Delta t} \,d\bx + \Psi^K + \sum_{\sigmat \in K}r_{\sigmat}  \\= \int_K\Big( \rho\dfrac{\delta^k \varepsilon}{\Delta t}  + p \xnabla\cdot\bvv\Big)\,\d\bx + \sum_{\sigmat \in K}r_{\sigmat}  = 0.
\end{multline}

The first term in \eqref{entropy:2} is a discrete analogue of \eqref{entropy:1}, therefore we can require 
\begin{equation*}
\int_K\Big( \rho\dfrac{\delta^k \varepsilon}{\Delta t}  + p \xnabla\cdot\bvv\Big)\,\d\bx \geq 0,
\end{equation*}
which yields another constraint on $r_{\sigmat}$:
\begin{equation}
\label{entropy_condition}
\sum_{\sigmat \in K}r_{\sigmat} \leq 0.
\end{equation}
We note that the derivation of the entropy condition for a general high-order scheme is slightly more tedious, however, it leads to exactly the same condition \eqref{entropy_condition} and is therefore not presented here.

Let us show that the entropy condition \eqref{entropy_condition} holds for the first order residual distribution scheme. From the conservation condition \eqref{correction:conservation} we have
\begin{equation*}
\label{correction:conservation-ord1}
r_{\sigmat} = \frac{1}{N_{\mathcal{E}}^K}\bigg( \int_{\partial K} \hat{p}_n \bvv\,d\sigma - \sum_{\sigmav \in K} \bvv_{\sigmav}\Phi_{\sigmav}^{K,\text{Rus}}  - \sum_{\sigmat \in K} \Psi_{\sigmat}^{K,\text{Rus}}  \bigg),
\end{equation*}
and therefore
\begin{equation}
\begin{split}
\label{correction:conservation-ord1-alpha-1}
\sum_{\sigmat \in K} r_{\sigmat} &= \int_{\partial K} \hat{p}_n \bvv\,d\sigma - \sum_{\sigmav \in K} \bvv_{\sigmav}\Phi_{\sigmav}^{K,\text{Rus}}  - \sum_{\sigmat \in K} \Psi_{\sigmat}^{K,\text{Rus}} \\ &=  \int_{\partial K} \hat{p}_n \bvv\,d\sigma - \sum_{\sigmav \in K} \bvv_{\sigmav} \Phi_{\sigmav}^K - \alpha_{K} \sum_{\sigmav \in K} \bvv_{\sigmav} (\bvv_{\sigmav}-\bar{\bvv})\\& \qquad \qquad \qquad \qquad- \sum_{\sigmat \in K} \Psi_{\sigmat}^K -
 \alpha_K \sum_{\sigmat \in K} (\varepsilon_{\sigmat}-\bar{\varepsilon}) \\ &=  \int_{\partial K} \hat{p}_n \bvv\,d\sigma - \int_K \bvv\cdot\xnabla p\,d\bx - \alpha_{K} \sum_{\sigmav \in K}  (\bvv_{\sigmav}-\bar{\bvv})^2 - \int_K p\xnabla\cdot\bvv\,d\bx  \\&= -\alpha_{K} \sum_{\sigmav \in K}  (\bvv_{\sigmav}-\bar{\bvv})^2 \leq 0,  
\end{split}\end{equation}
where we have taken into account that 
\begin{equation*}
\sum_{\sigmav \in K} \bvv_{\sigmav} (\bvv_{\sigmav}-\bar{\bvv}) = \sum_{\sigmav \in K}  (\bvv_{\sigmav}-\bar{\bvv})^2, \quad \text{and}\quad \sum_{\sigmat \in K} (\varepsilon_{\sigmat}-\bar{\varepsilon}) = 0.
\end{equation*}
Therefore, the entropy condition \eqref{entropy_condition} is satisfied with any $\alpha_{K}  \geq 0$.


\section{Optimization of artificial viscosity}
\label{sec:viscosity}

The artificial viscosity coefficient $\alpha_{K}$ present in the first-order Rusanov residual \eqref{Eq:residual-alpha-velo} plays a crucial role in ensuring the stability of high order staggered FEM approximation, and it defines the amount of entropy dissipation as follows from relation \eqref{correction:conservation-ord1-alpha-1}. Therefore, on one hand, excessively large values of this coefficient will stabilize the method but, on the other hand, will lead to excessive numerical diffusion of the solution features. This might be more critical for problems involving vortical flows since the numerical dissipation will deteriorate the resolution of the high order finite elements and prevent the development of physically correct vortical structures in the numerical solution. 

Therefore, in order to reduce the amount of numerical dissipation we adopt a MARS (Multidirectional Approximate Riemann Solution) technique proposed in \cite{MARS2014,MARS2017} for the construction of the artificial viscosity term. The idea of MARS approach is based on considering a multidirectional Riemann problem at the nodes of the mesh element and using the solution of this problem for the approximation of forces acting on every node. 

The original artificial viscosity term typically used in residual distribution schemes has the form 
\begin{equation}
\label{eq:LxF_visc}
\sigma_a^{Rus} = \alpha_{K}(\bvv_{\sigmav}-\bar{\bvv}),
\end{equation}
where $\bar{\bvv} = \dfrac{1}{N_{\mathcal{V}}^K} \sum_{{\sigmav}\in K} \bvv_{\sigmav}$ and $\alpha_{K}$ is an estimate of the largest eigenvalue of the system and is defined by the shock impedance $\rho U$ multiplied by some length scale of the element, regardless of the direction of the flow and the number of DOF $\sigmav \in K$. MARS approach allows to put a sensor on the artificial viscosity term so that different amount is added at different DOFs inside the cell $K$.

The MARS artificial viscosity term will be defined as 
\begin{equation}
\label{eq:MARS_visc}
\sigma_a^{MARS} = \alpha_{K}^i (\bvv_{\sigmav}-\tilde{\bvv}),
\end{equation}
where $\alpha_K^i = \rho U |\be_i\cdot\bn_i|$ and $\tilde{\bvv} = \sum_{\sigmav\in K} \alpha_K^i \bvv_{\sigmav }/ \sum_{\sigmav\in K} \alpha_K^i$. 
Using $\sigma_a^{MARS} $, the entropy balance \eqref{correction:conservation-ord1-alpha-1} becomes
$$\sum_{\sigmat \in K} r_{\sigmat}  = - \sum_{\sigmav \in K} \alpha_{K}^i (\bvv_{\sigmav}-\tilde{\bvv})^2 \leq 0,$$
and one has to define $\alpha_{K}^i$ such that
$$\sum_{\sigmav \in K} \alpha_{K}^i (\bvv_{\sigmav}-\tilde{\bvv})^2 \leq \alpha_K\sum_{\sigmav \in K}  (\bvv_{\sigmav}-\bar{\bvv})^2 .$$

The vector $\be_i$ is a unit vector which approximates the direction of the shock and the vector $\bn_i$ is defined by $\bn_i = \int_K \nabla\varphi_i\,d\bx$. Therefore, the maximal numerical viscosity will be applied in the direction of the shock. In \cite{MARS2017}, the shock direction is chosen as $\be_i = \frac{\bvv_{\sigmav}-\bar{\bvv}}{\|\bvv_{\sigmav}-\bar{\bvv}\|}$,  however, other choices are possible. For example, we set $\be_i = \frac{\bvv_{\sigmav}}{\|\bvv_{\sigmav}\|}$ for the Sedov and Noh problems since these problems have radial symmetry and therefore the direction of maximal compression is aligned with the velocity. Finally, the impedance $\rho U$ is calculated as $\rho U = \rho (c + s||\Delta \bu||)$, $\Delta\bu = \bvv_{\sigmav}-\bar{\bvv}$ and $s = \frac{\gamma+1}{2}$. 

At this point we would like to highlight the differences between the artificial viscosity approach in \cite{rieben,MARS2014,MARS2017} and the one proposed here. The main difference consists in ways to achieve high order of accuracy: thus, in \cite{MARS2014,MARS2017}, a modification to \eqref{eq:MARS_visc} would be needed in order to transition from first to higher orders. Contrary to that, the philosophy of RD method is such that only first order viscosity is sufficient in order to achieve high order. This is because the first order Rusanov residuals are used to calculate the distribution coefficients $\beta_i^K$ according to \eqref{Psi2}, but the high order residuals are defined as distributions of the total residual as stated by \eqref{RD:Phi_unsteady}, \eqref{def_phiK} and hence the order of approximation can be preserved.

The viscosity term in \eqref{Eq:residual-alpha-eps} is modified in a similar way.

\section{Numerical results}
\label{sec:numerical_results}

In this section, we apply the multidimensional SGH RD scheme to several well-known test problems in Lagrangian hydrodynamics to assess its robustness and accuracy. We perform the simulations using the second-order SGH RD scheme which is based on quadratic Bernstein shape functions for the approximation of kinematic variables and piecewise-linear shape functions for the thermodynamic variables, and the second-order timestepping algorithm described in Section~\ref{sec:timestepping}. Finally, unless stated otherwise, we use the MARS artificial viscosity from Section~\ref{sec:viscosity}.

\subsection{Taylor-Green vortex}

The Taylor-Green vortex problem is typically used for the assessment of the accuracy of the Lagrangian solvers \cite{rieben}. The purpose of this test case is to verify the
ability of the fully discrete methods to obtain high-order convergence in time and space
on a moving mesh with nontrivial deformation for the case of a smooth problem. Here we consider a simple, steady state solution to the 2D incompressible, inviscid Navier–Stokes equations, given by the initial conditions
\[ \bu = \{ \sin(\pi x) \cos(\pi y), −\cos(\pi x) \sin(\pi y)\}, \quad p = \dfrac{\rho}{4}\big(\cos(2\pi x) + \cos(2\pi y)\big) + 1. \]

We can extend this incompressible solution to the compressible case with an ideal
gas equation of state and constant adiabatic index $\gamma = 5/3$ by using a manufactured
solution, meaning that we assume these initial conditions are steady state solutions
to the Euler equations, then we solve for the resulting source terms and use these
to drive the time-dependent simulation. The flow is incompressible $(\nabla\cdot \bu = 0)$ so
the density field is constant in space and time and we use $\rho \equiv 1$. It is easy to check
that $\rho \diff{\bu}{t} = ˆ'\nabla p$ so the external body force is zero. In the energy equation, using
$\varepsilon = p/((\gamma - 1)\rho)$, we compute
\[ \varepsilon_{\mathrm src} = \rho \diff{\varepsilon}{t} + p\nabla\cdot\bu = \diff{\varepsilon}{t} = \dfrac{3\pi}{8}(\cos(3\pi x) \cos(\pi y) − \cos(\pi x) \cos(3\pi y)).\]

This procedure allows us to run the time-dependent problem to some point in time,
then perform normed error analysis on the final computational mesh using the exact
solutions for $\bu$ and $p$. The computational domain is a unit box with wall boundary
conditions on all surfaces $(\bu \cdot \bn = 0)$. Note that for this manufactured solution all
fields are steady state, i.e., they are independent of time; however, they do vary along
particle trajectories and with respect to the computational mesh as it moves. We run
the problem until $T = 1.25$. Since this problem is smooth we run it without any artificial
viscosity (i.e. we set the artificial viscosity of the Rusanov or MARS scheme to $0$) and do normed error analysis on the solution at the final time and compute
convergence rates using a variety of high-order methods.

\begin{figure}[ht]
 \begin{center}	
	\begin{tabular}{ccc}
	\hspace{-0.5cm}\includegraphics[width=0.29\textwidth]{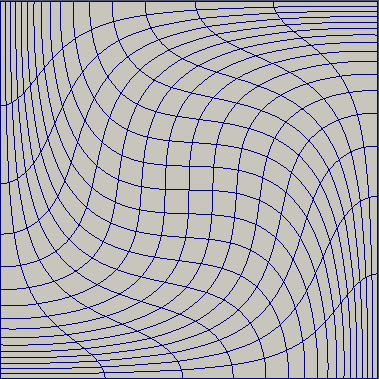} & \hspace{-0.4cm}\includegraphics[width=0.29\textwidth]{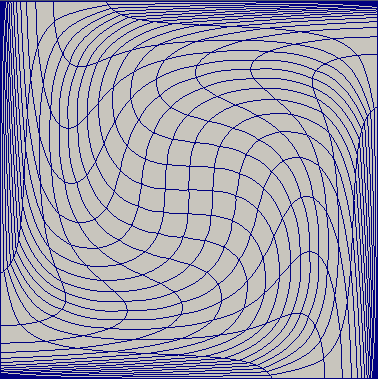} & \hspace{-0.4cm}\includegraphics[width=0.29\textwidth]{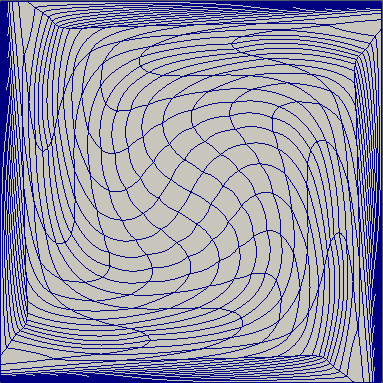} \\
    \hspace{-0.7cm}\includegraphics[width=0.35\textwidth]{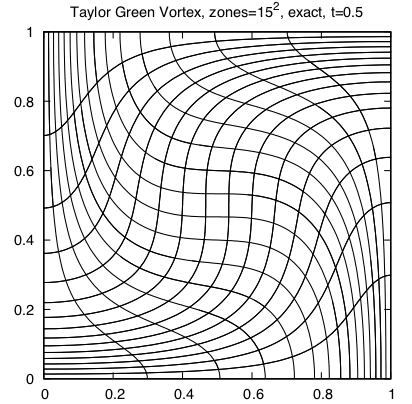} & \hspace{-0.6cm}\includegraphics[width=0.35\textwidth]{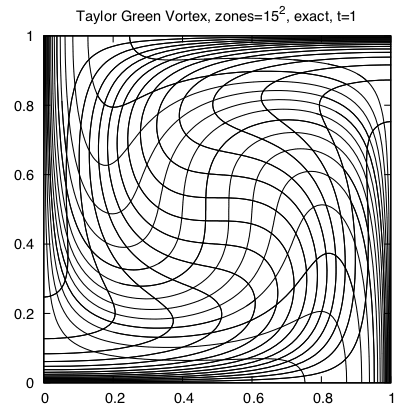} & \hspace{-0.6cm}\includegraphics[width=0.35\textwidth]{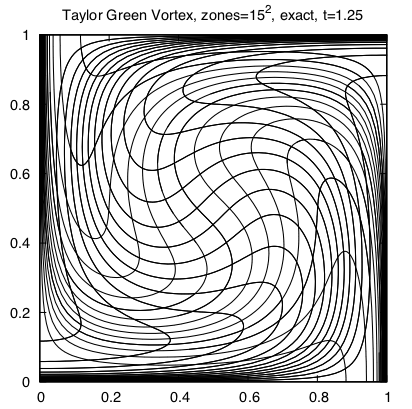} \\ 
    $t=0.5$ & $t=1.0$ & $t=1.25$
	\end{tabular}
 \end{center}	
 \caption{Taylor-Green vortex, numerical (top) and exact (bottom)}
 \label{Fig:TGV}
\end{figure}

\begin{figure}[ht]
	\begin{center}
		\begin{tabular}{cc}
			\includegraphics[width=0.45\textwidth]{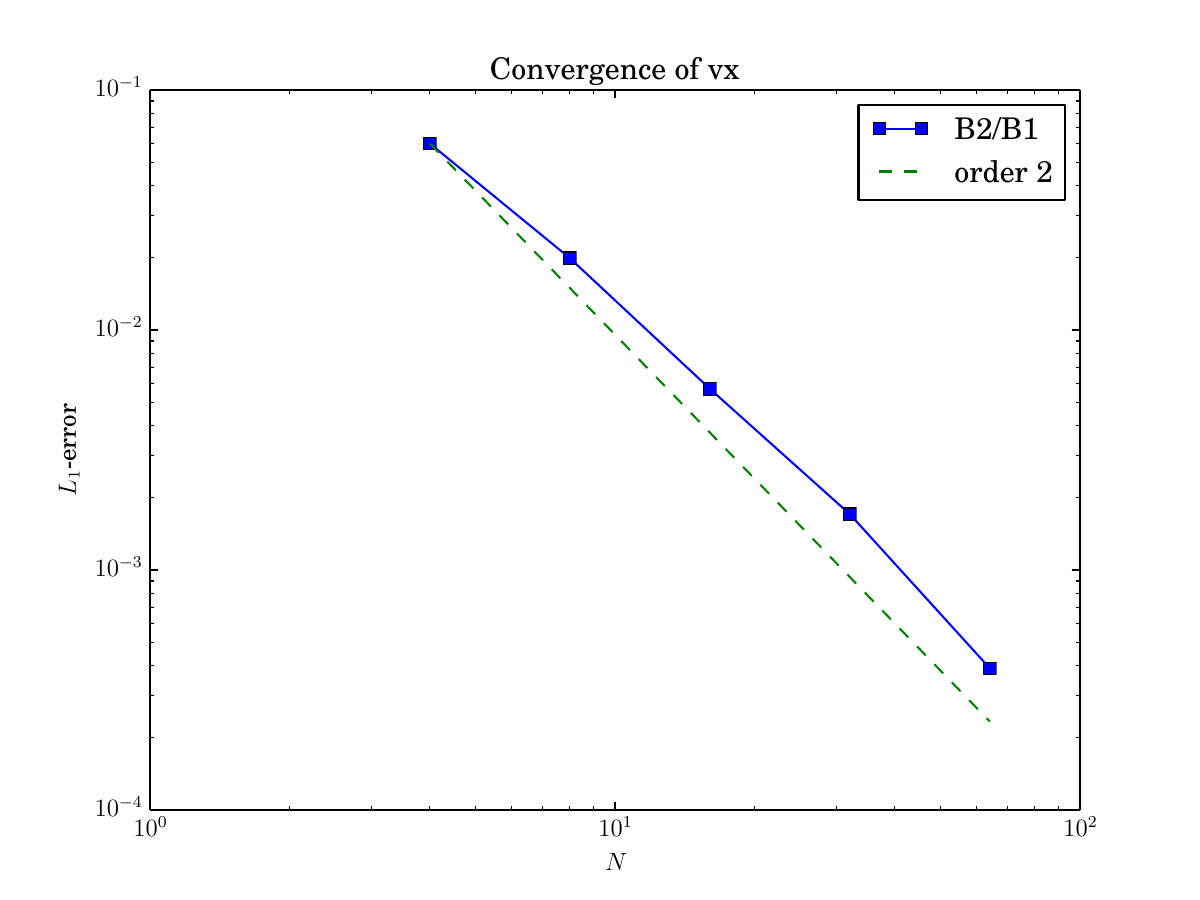} &  \includegraphics[width=0.45\textwidth]{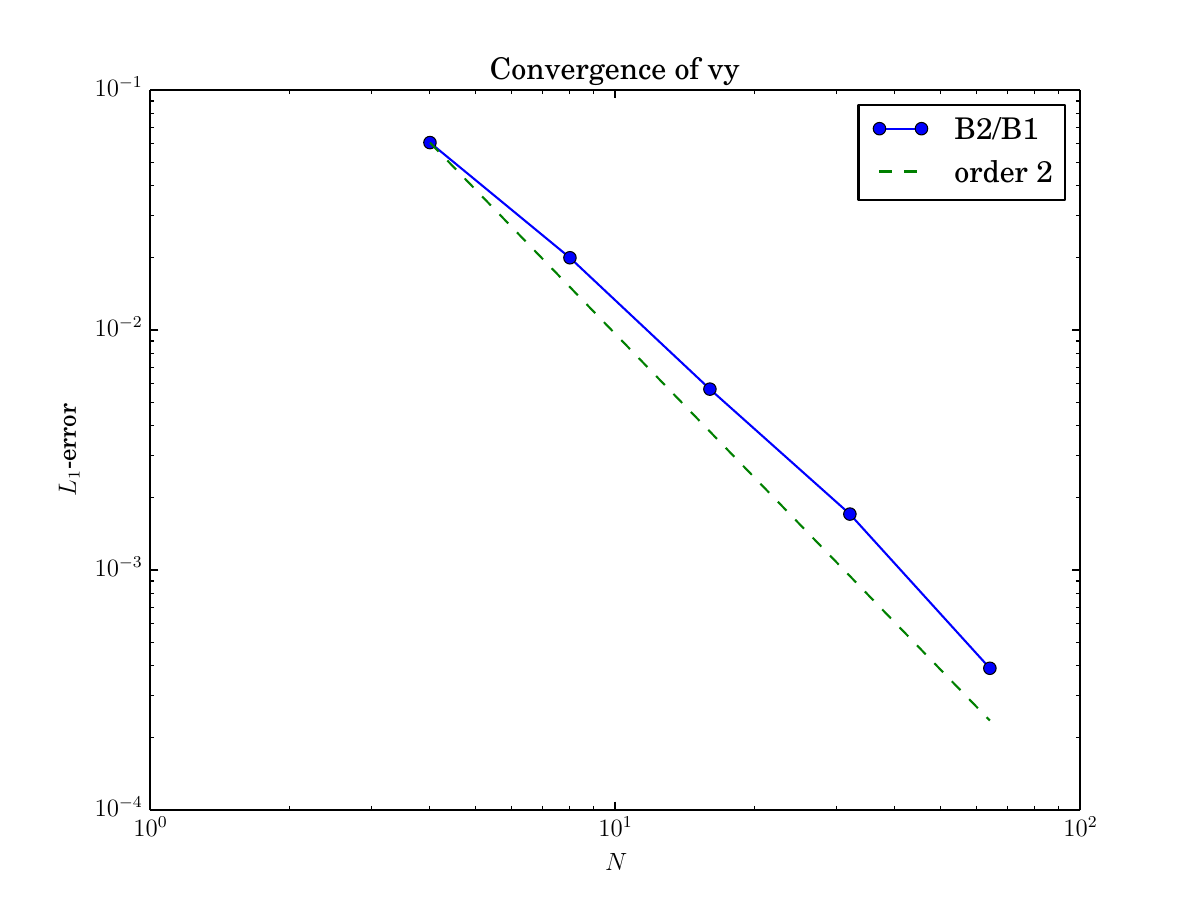}  \\
			Convergence of $v_x$ & Convergence of $v_y$ 
		\end{tabular}
	\end{center}
    \caption{Taylor-Green vortex, convergence at $t=0.5$}
    \label{Fig:TGV_convergence}
\end{figure}

In Fig.~\ref{Fig:TGV} we show plots of the curvilinear mesh at times $t = 0.5$, $t = 1.0$ and $t = 1.25$, and we compare the numerical (upper row) solution with the exact one (lower row). In Fig.~\ref{Fig:TGV_convergence}, we plot the errors of the velocity components in $L_1$ norm vs the mesh resolution.

\subsection{Gresho vortex}

The Gresho problem \cite{Gresho1990,LiskaWendroff2003} is a rotating vortex problem
independent of time. Angular velocity depends only on radius, and the centrifugal
force is balanced by the pressure gradient
\begin{equation*}
u_\phi(r) =
\begin{cases}
5r, \\
2 - 5r, \\
0 \\
\end{cases}
p(r) =
\begin{cases}
5 + \frac{25}{2} r^2, \quad &0 \leq r < 0.2, \\
9 - 4 \ln 0.2 + \frac{25}{2} r^2 - 20r + 4\ln r, \quad &0.2 \leq r < 0.4, \\
3 + 4\ln 2, \quad &0.4 \leq r.
\end{cases}
\end{equation*}
The radial velocity is $0$ and the density is $1$ everywhere. 

\begin{figure}[ht]
	\begin{center}
	\begin{tabular}{cc}
		\includegraphics[width=0.46\textwidth]{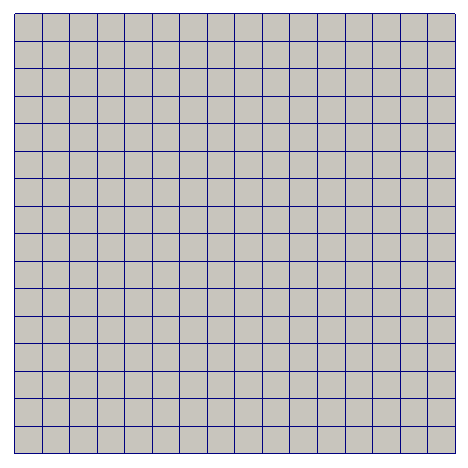} & \includegraphics[width=0.45\textwidth]{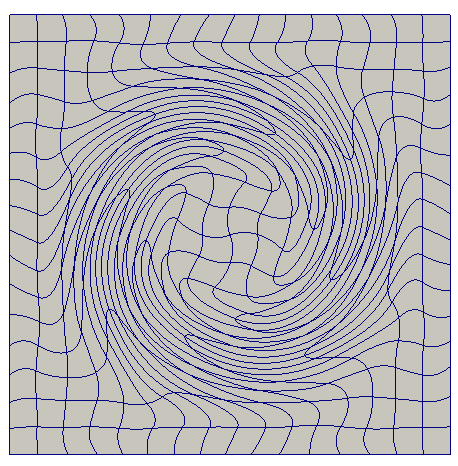}
	\end{tabular}
	\end{center}
	\caption{Gresho vortex, numerical solution, mesh at $t=0$ and $t=0.65$}
	\label{Fig:Gresho}
\end{figure}

Gresho problem is an interesting validation test case to assess the robustness and the accuracy of a Lagrangian scheme. The vorticity leads to a strong mesh deformation which can cause some problems such as negative Jacobian determinants or negative densities. We compute the solution to this problem on a rectangle $[0, 1]\times[0, 1]$ until time $T = 0.65$ on a $16\times16$ grid using the second-order SGH RD scheme. The initial and final grids are shown in Fig.~\ref{Fig:Gresho}. We observe that our scheme is able to evolve the vortex robustly until the mesh becomes strongly tangled.

\subsection{Sedov problem}

Next, we consider the Sedov problem for a point-blast in a uniform medium. An exact solution based on self-similarity arguments is available, see for instance \cite{kamm,Sedov}. This test problem provides a good assessment of the robustness of numerical schemes for strong shocks as well as the ability of the scheme to preserve cylindrical symmetry. 

%
%

\begin{figure}[ht]
	\begin{center}
		\begin{tabular}{cc}
			\includegraphics[width=0.5\textwidth]{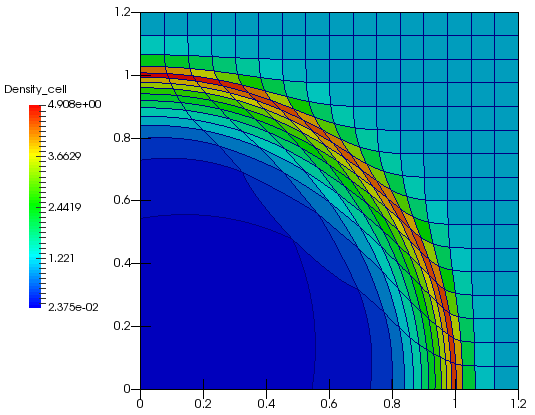} & \hspace{-0.5cm}\includegraphics[width=0.45\textwidth]{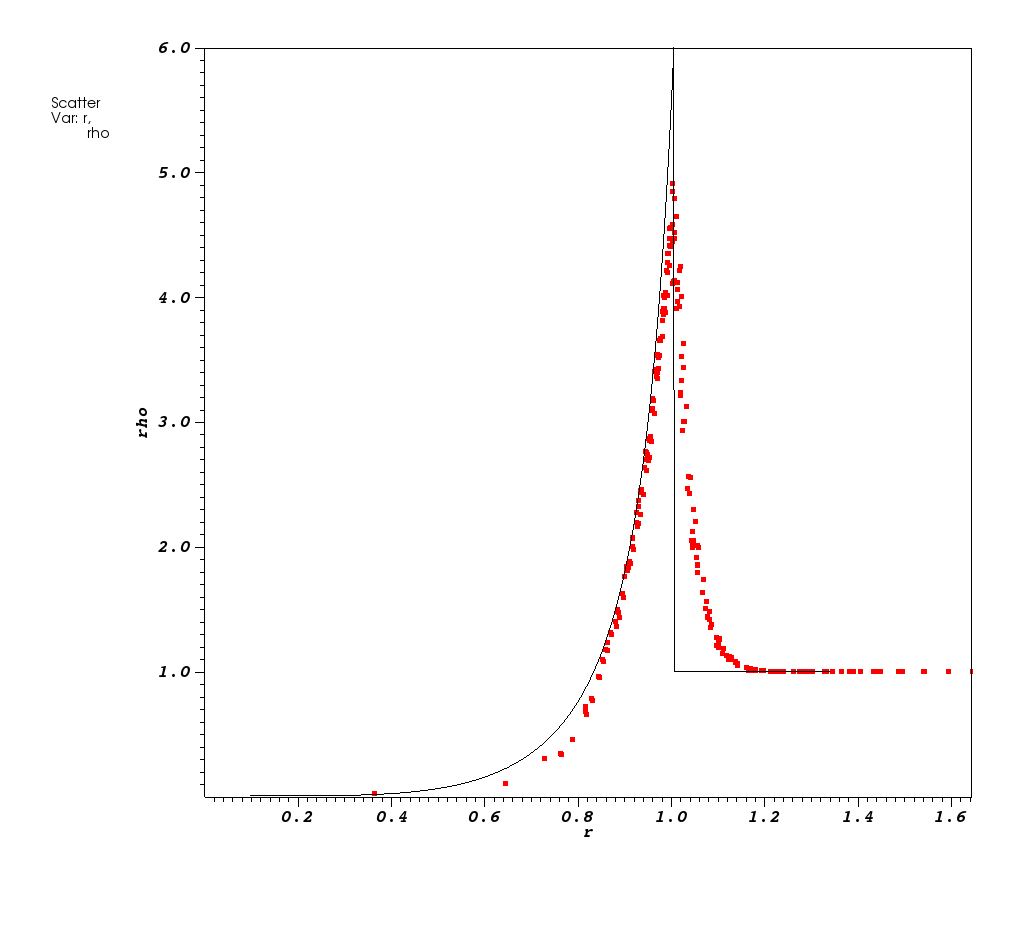}
		\end{tabular}
	\end{center}
	\caption{Sedov problem, $16\times16$ mesh, density at $t=1.0$}
	\label{Fig:Sedov_16x16_MARS}
\end{figure}

\begin{figure}[ht]
	\begin{center}
		\begin{tabular}{cc}
			\includegraphics[width=0.53\textwidth]{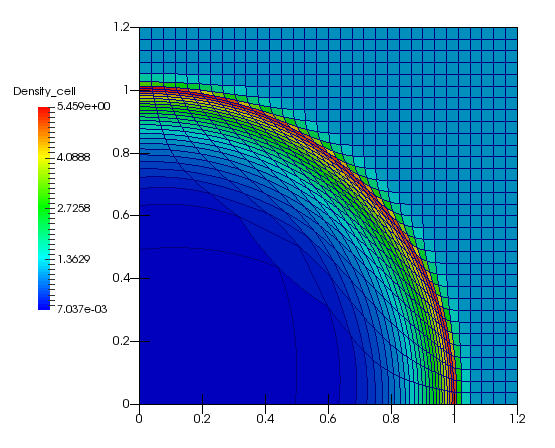} & \hspace{-0.6cm}\includegraphics[width=0.45\textwidth]{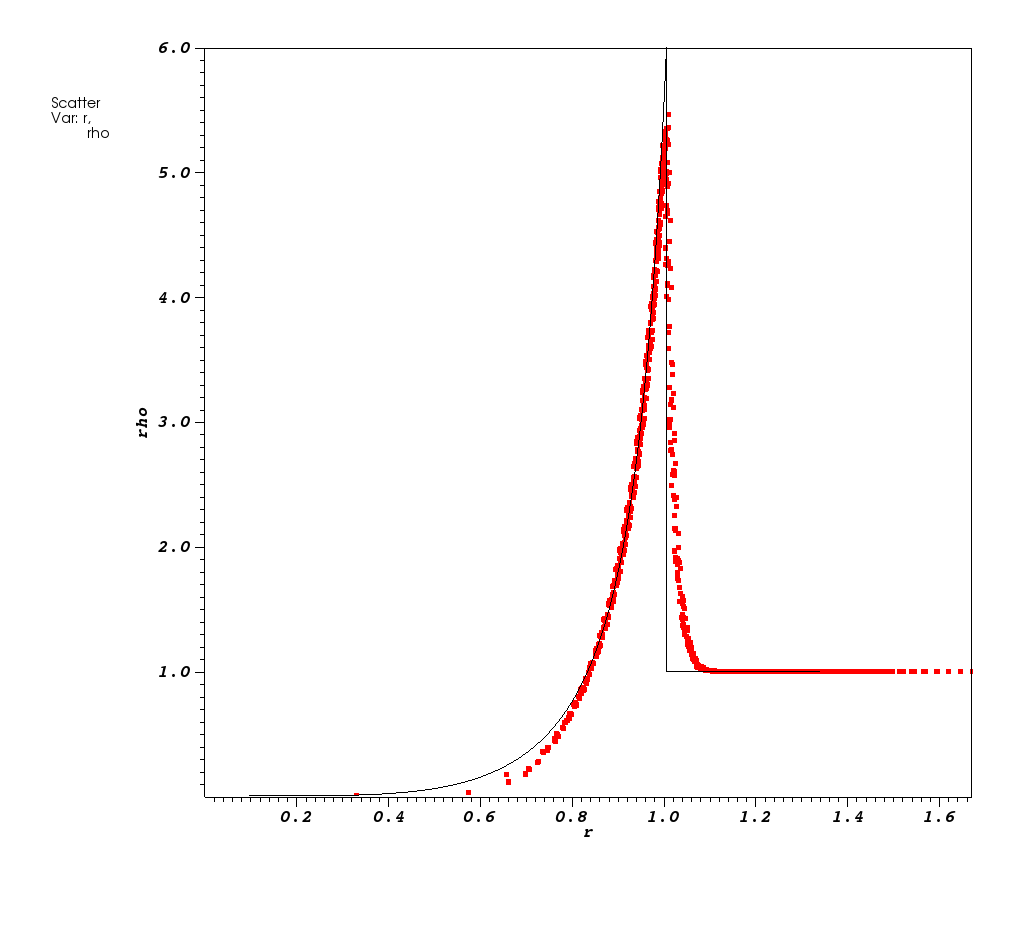}
		\end{tabular}
	\end{center}
	\caption{Sedov problem, $32\times32$ mesh, density at $t=1.0$}
	\label{Fig:Sedov_32x32_MARS}
\end{figure}

The Sedov problem consists of an ideal gas with $\gamma = 1.4$ and a delta source of internal energy imposed at the origin such that the total energy is equal to $1$. The initial data is $\rho_0 = 1$, $u_0 = v_0 = 0$, $p_0 = 10^{-6}$. At the origin, the pressure is set to 
\[
p_0 = \dfrac{(\gamma-1)\rho_0\varepsilon_S}{V_{or}}, 
\]
where $V_{or}$ is the volume of the cell containing the origin and $\varepsilon_S = 0.244816$ as suggested in \cite{kamm}. The solution consists of a diverging infinite strength shock wave whose front is located at radius $r=1$ at time $T=1$, with a peak density reaching $6$. 
We first run the Sedov problem with the second-order SGH RD scheme with MARS viscosity on a $16\times16$ Cartesian grid in the domain $[0, 1.2]\times[0, 1.2]$. The results are shown in Fig.~\ref{Fig:Sedov_16x16_MARS}. At the end of the computation, the shock wave front is correctly located and is symmetric. The density peak reaches $4.908$ with MARS viscosity, which we consider to be a very good approximation on such coarse grid. Next, we run the same test case on a finer grid consisting of $32\times32$ cells, the corresponding results are presented in Fig.~\ref{Fig:Sedov_32x32_MARS}. Note that the peak density value $5.459$ becomes closer to the exact value by mesh refinement.


These results demonstrate the robustness and the accuracy of our  scheme.

\subsection{Noh problem}

\begin{figure}[ht]
	\begin{center}
		\includegraphics[scale=0.2]{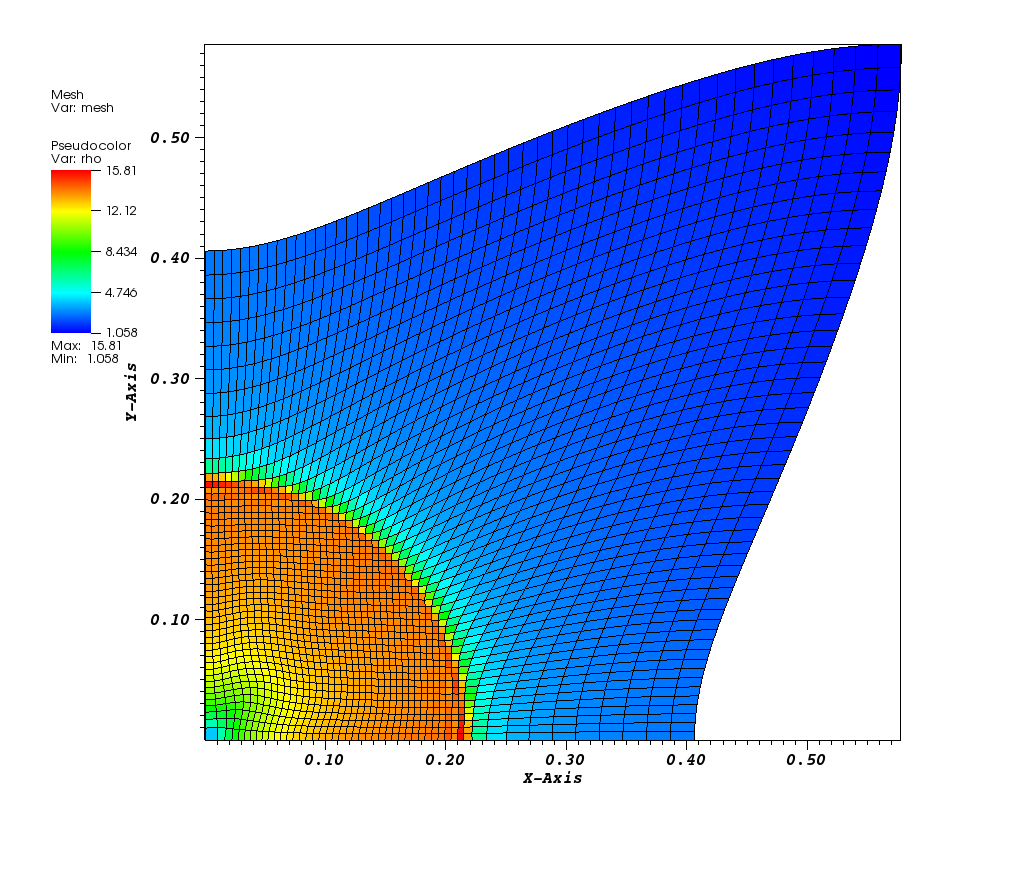} \\
		\begin{tabular}{cc}
			\includegraphics[scale=0.17]{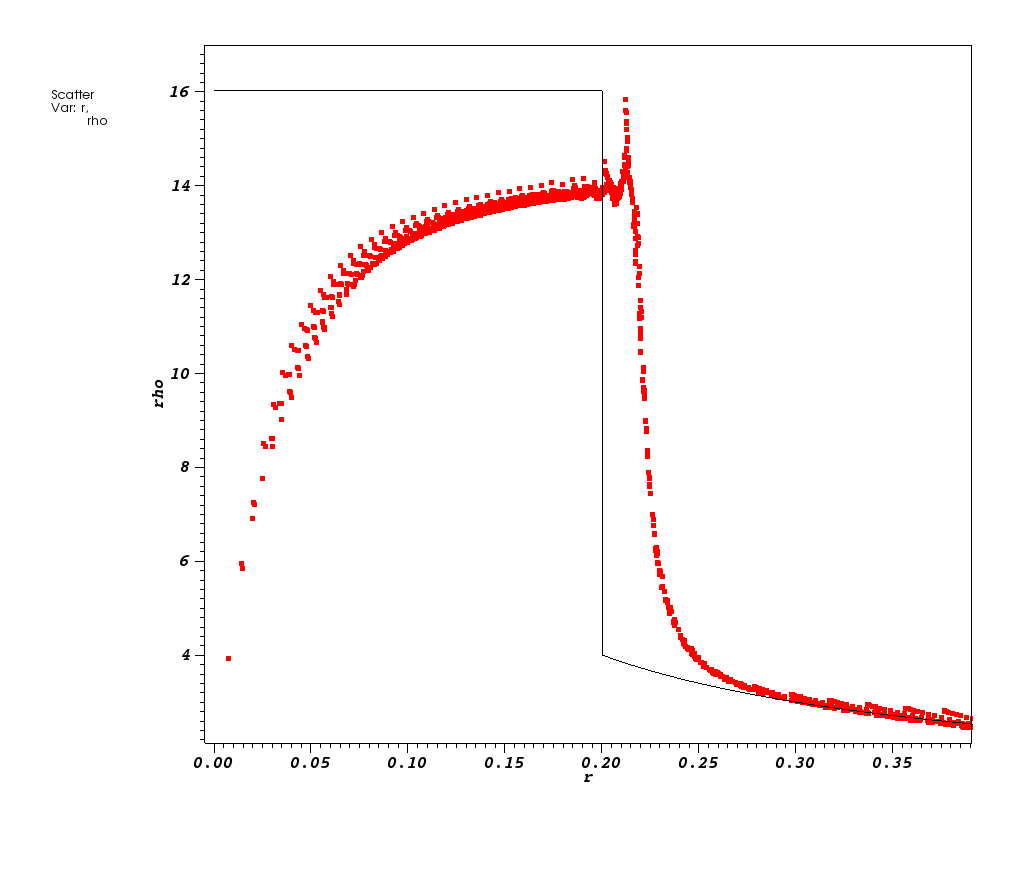} & \hspace{-0.5cm}\includegraphics[scale=0.17]{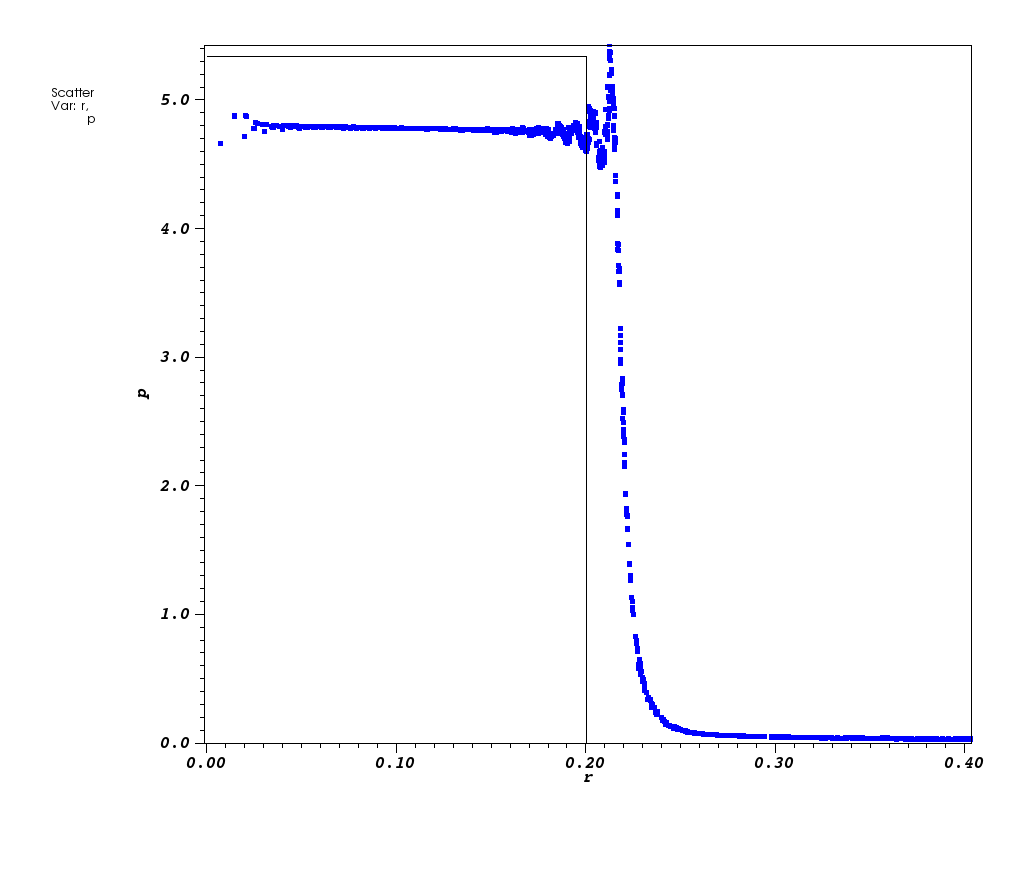} 
		\end{tabular}	
	\end{center}
	\caption{Noh problem, $50\times50$ mesh, density contour (top), density scatter (bottom left) and pressure scatter (bottow right) at $t=0.6$}
	\label{Fig:Noh}
\end{figure}

The Noh problem \cite{noh} is a well known test case used to validate Lagrangian schemes in the regime of infinite strength shock wave. The problem consists of an ideal gas with $\gamma = 5/3$, initial density $\rho_0 = 1$,
and initial energy $\varepsilon_0 = 0$. The value of each velocity degree of freedom is initialized to
a radial vector pointing toward the origin, $\bu = \mathbf{r}/\|\mathbf{r}\|$. The initial velocity generates
a stagnation shock wave that propagates radially outward and produces a peak post-shock
density of $\rho = 16$. The initial computational domain is defined by $[0, 1]\times[0, 1]$. We run the Noh problem on a $50\times50$ Cartesian grid using the SGH RD scheme with MARS artificial viscosity. This configuration leads to a severe test case since the mesh is not aligned with the flow. In Fig.~\ref{Fig:Noh} we show plots of the density field at the final time of $t = 0.6$ as well as scatter plots of density versus
radius.
We note that we have a very smooth and cylindrical solution, and that the shock is located at a circle whose radius is approximately $0.2$. We see that the numerical solution is very close to the one-dimensional analytical solution, and the numerical post-shock density is not too far from the analytical value. This shows the ability of our scheme to preserve the radial symmetry of the flow.


\subsection{Triple point problem}

\begin{figure}[ht] 
	\centering
	\includegraphics[scale=0.25]{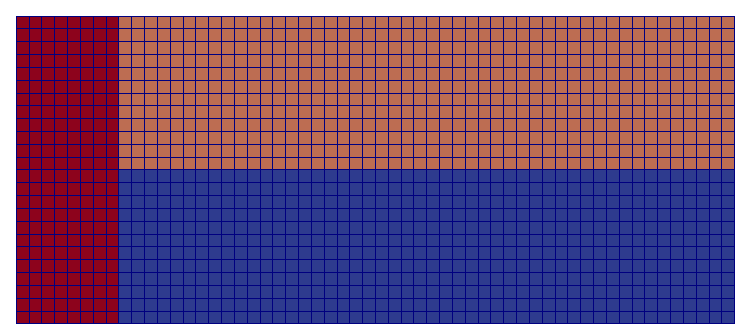}
	\caption{\label{TriplePointInitial}  The initial mesh for the triple point problem is shown above.  The initial conditions are as follows: left region (red) is $\rho = 1.0 \, \& \, e = 2.0$, top-right region (pink) is $\rho = 0.125 \, \& \, e=1.6$, and the bottom-right region (purple) is $\rho = 1.0 \, \& \, e=0.25$. }  
\end{figure}

\begin{figure}[ht]
	\begin{center}	
		\begin{tabular}{cc}
		  \hspace{-0.5cm}\includegraphics[width=0.52\textwidth]{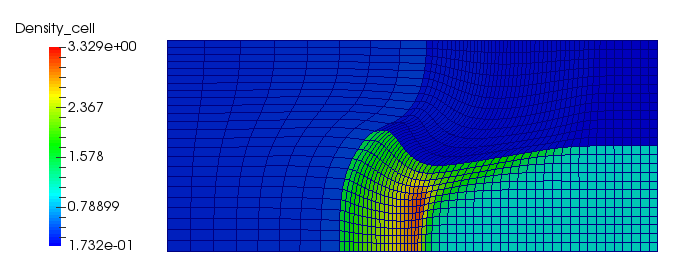} &  \hspace{-0.5cm}\includegraphics[width=0.52\textwidth]{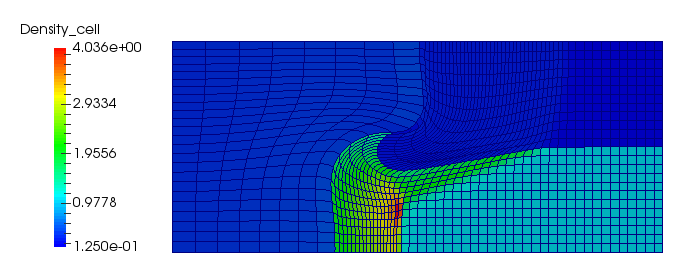} \\
		  $\sigma_a^{Rus}$, $1$st order & $\sigma_a^{MARS}$, $1$st order \\
		   \hspace{-0.5cm}\includegraphics[width=0.52\textwidth]{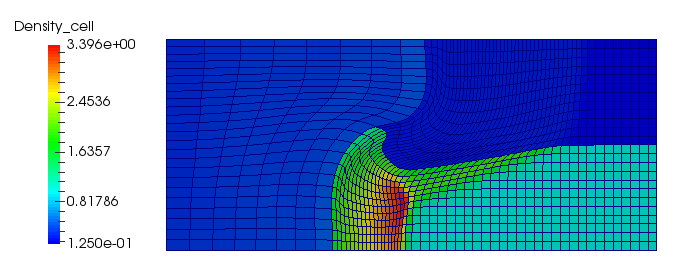} &  \hspace{-0.5cm}\includegraphics[width=0.52\textwidth]{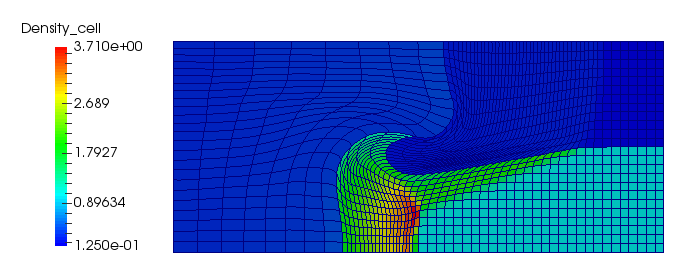} \\
		  $\sigma_a^{Rus}$, $2$nd order & $\sigma_a^{MARS}$, $2$nd order 
		\end{tabular}
	\end{center}	
    \caption{Triple point problem, $24\times56$ mesh, density at $t=3.0$}
    \label{Fig:TriplePoint}
\end{figure}

The triple point problem is used to assess the robustness of a Lagrangian method on a problem that has significant vorticity \cite{Kucharik}. The initial conditions are three regions of a gamma-law gas, where each region has different initial conditions. One region has a high pressure that drives a shock through two connected regions and a vortex develops at the triple point where three regions connect.  In this study, every region uses a gamma of $1.4$.  Fig. \ref{TriplePointInitial} shows the initial conditions, and Fig. \ref{Fig:TriplePoint} shows the results at $3.0 \mu s$.  The mesh remained stable despite large deformation and calculation will continue to run well beyond $5 \mu s$.  The triple point problem illustrates that SGH RD method can be used on problems with significant mesh deformation. To demonstrate the effect of viscosity optimization, we run this problem both using $\sigma_a^{Rus}$ and $\sigma_a^{MARS}$ from Section~\ref{sec:viscosity}.

\section{Conclusions}
\label{sec:conclusions}

In this paper we have extended a Residual Distribution (RD) scheme for the Lagrangian hydrodynamics to multiple space dimensions and proposed an efficient way to construct artificial viscosity terms. We have developed an efficient mass matrix diagonalization algorithm which relies on the modification of the time-stepping scheme and gives rise to an explicit high order accurate scheme. 
The two-dimensional numerical tests considered in this paper show the robustness of the method for problems involving very strong shock waves or vortical flows.
	
Further research includes the extension of the present method to higher order in space and time. We also plan to extend the method to solids.

\section*{Acknowledgments} This research has been supported by Laboratory Directed Research and Development (LDRD) program at Los Alamos National Laboratory. The Los Alamos unlimited release number is LA-UR-18-30342. ST and RA also thank the partial financial support of the Swiss National Science Foundation (SNSF) via the grant \# 200021\_153604 ("High fidelity simulation for compressible material").  We acknowledge several discussions with D. Kuzmin (TU Dortmund, Germany) about the domain invariance properties of the Rusanov scheme.

\bibliographystyle{siamplain}
\bibliography{biblio}
\appendix
\input{conservation}
\end{document}

%% file: conservation.tex
\section{Justification of \eqref{Eq:conserv-discr}}\label{appendix:conservation}
The time stepping method writes as \eqref{eq:RD_scheme} (if we forget about the mesh movement):
	\begin{equation}
	\label{Eq:RD_scheme-velo:1}
	C_{\sigmav}^{\mathcal{V}}\dfrac{\bvv_{\sigmav}^{(k+1)} - \bvv_{\sigmav}^{(k)}}{\Delta t}  + \sum_{K \ni \sigmav}
	\Phi_{{\sigmav},ts}^{K,L} = 0, 
	\end{equation}
	\begin{equation}
	\label{Eq:RD_scheme-eps:1}
	C_{ \sigmat}^{\mathcal{E}}\dfrac{\varepsilon_{\sigmat}^{(k+1)}-\varepsilon_{\sigmat}^{(k)}}{\Delta t}  + \sum_{K \ni \sigmat }
	\Psi_{{\sigmat},ts}^{K,c} = 0,
	\end{equation}
In \cite{AbgrallRoe2001} was given a Lax-Wendroff type result that guaranties, under standard stability requirements, that provided a conservation relation at the level of the element is satisfied, if the sequence of approximate solutions converges to a limit in some $L^p$ space, is a weak solution of the discretized hyperbolic system. The proof uses the fact that in each mesh element, the solution is approximated by a function from a finite dimensional approximation space. It can easily be extended to unsteady problems, and in that case, the total residual has a time increment component, and a spatial term.  If we were integrating the conservative system
$$\dpar{\mathbf{u}}{t}+\text{ div }\mathbf{f}(\mathbf{u})=0,$$
the total residual, i.e. a consistent approximation of
$$\int_{t_n}^{t_{n+1}} \int_K\bigg (\dpar{\mathbf{u}}{t}+\text{ div }\mathbf{f}(\mathbf{u})\bigg ) \;dx\; dt,$$
would be the sum of a time increment
$$\Phi_t =\int_K \big (\mathbf{u}^{n+1}-\mathbf{u}^n\big ) \; dx  \approx \int_{t_n}^{t_{n+1}} \int_K\dpar{\mathbf{u}}{t}\;dx\; dt$$
and a space term
$$\Phi_x\approx \int_{t_n}^{t_{n+1}} \int_K \text{ div }\mathbf{f}(\mathbf{u})\;dx\; dt .$$
It turns out that we can be quite flexible in the approximation of  $\Phi_t $, provided it remains the difference between a term evaluated at $t_{n+1}$ and a term evaluated at $t_n$, but we must be very careful in the definition of the spatial term: all this is very similar to what happens for the classical Lax-Wendroff theorem. In addition, the time stepping we are using is set in such a way that the time increment has the form
$$\int_K \big (\mathbf{u}^{l+1}-\mathbf{u}^l\big ) \; dx$$ with $l=0,1$, and the space increment contains also a time increment between the iterations $0$ and $l$. However this does not change anything fundamentally.

In our case, despites  the staggered nature of the grid, we can define a specific energy on the element $K$: we have a specific internal energy \emph{function} $\varepsilon_h$ defined on $K$, and a velocity field $\bvv_h$ also defined on $K$: hence we define $e_h=\varepsilon_h+\tfrac{1}{2}(\bvv^h)^2$, and the increment in time of the energy would be
$$\int_K \rho_h^{n+1} \big ( \varepsilon_h^{n+1}+\tfrac{1}{2}(\bvv_h^{n+1})^2\big ) \; dx-\int_K \rho_h^{n} \big ( \varepsilon_h^{n}+\tfrac{1}{2}(\bvv_h^{n})^2\big ) \; dx.$$
However doing this, it is very complicated to estimate this quantity. Since we can be flexible, if we keep the incremental structure, we can approximate the variation of energy between two cycles by
\begin{equation}\label{appendix:1}
\begin{split}
\delta e&=\sum_{i_{\mathcal{E}}\in K}C_{i_\varepsilon}^{\mathcal{E},K} (\varepsilon_{i_{\mathcal{E}}}^{l+1}-\varepsilon_{i_{\mathcal{E}}}^{l}\big )+
\sum_{i_\mathcal{V}\in K} C_{i_\mathcal{V}}^{\mathcal{V},K}\;\frac{1}{2} \bigg (\big ( \bvv_{i_\mathcal{V}}^{l+1}\big )^2-\big ( \bvv_{i_\mathcal{V}}^{l}
\end{split}\big )^2\bigg )\end{equation}
This relation can be rewritten as
$$\delta e=\sum_{i_{\mathcal{E}}\in K}C_{i_\varepsilon}^{\mathcal{E},K} (\varepsilon_{i_{\mathcal{E}}}^{l+1}-\varepsilon_{i_{\mathcal{E}}}^{l}\big )+
\sum_{i_\mathcal{V}\in K} C_{i_\mathcal{V}}^{\mathcal{V},K} \;\tilde{\bvv}_{i_\mathcal{V}}\cdot \bigg ( \bvv_{i_\mathcal{V}}^{l+1}-\bvv_{i_\mathcal{V}}^{l}\bigg )$$
where
$$\tilde{\bvv}_{i_\mathcal{V}}=\frac{1}{2}\big ( \bvv_{i_\mathcal{V}}^{l+1}+ \bvv_{i_\mathcal{V}}^{l}\big ).$$

Now the idea is to multiply \eqref{Eq:RD_scheme-velo:1} by $\tilde{\bvv}_{i_\mathcal{V}}$, and add the relation \eqref{Eq:RD_scheme-eps:1}. Hence, we get:
\begin{equation*}
\begin{split}
\mathcal{D}&=\sum_{i_{\mathcal{E}}\in K}C_{i_\varepsilon}^{\mathcal{E},K} (\varepsilon_{i_{\mathcal{E}}}^{l+1}-\varepsilon_{i_{\mathcal{E}}}^{l}\big )+
\sum_{i_\mathcal{V}\in K} C_{i_\varepsilon}^{\mathcal{V},K} \tilde{\bvv}_{i_\mathcal{V}}\cdot \bigg ( \bvv_{i_\mathcal{V}}^{l+1}-\bvv_{i_\mathcal{V}}^{l}\bigg )\\
&+
\Delta t \bigg ( \sum_{i_{\mathcal{E}}\in K}\Psi_{{\sigmat},ts}^{K,c}+\sum_{i_\mathcal{V}\in K}  \tilde{\bvv}_{i_\mathcal{V}}\cdot \Phi_{{\sigmav},ts}^{K,L}\bigg )
\end{split}
\end{equation*}
What we ask is that this quantity be \emph{equal} to a suitable approximation of 
$$\int_K \big ({e}^{l+1}-{e}^l\big ) \; dx +\int_K \big (e^l-e^0\big )\; dx +\Delta t \int_{\partial K}\dfrac{p^{(l)} \bvv^{(l)} +p^{(0)}\bvv^{(0)}}{2}\cdot \bn \; d\sigma.$$
Of course we approximate 
$\int_K \big ({e}^{l+1}-{e}^l\big ) \; dx$ as \eqref{appendix:1}, and $\int_K \big (e^l-e^0\big )\; dx$ similarly for simplicity, so that we end up with the relation
\begin{equation*}\label{appedix:conservationEnergy}
\begin{split}
\sum_{i_{\mathcal{E}}\in K}\Psi_{{\sigmat},ts}^{K,c}+\sum_{i_\mathcal{V}\in K}  \tilde{\bvv}_{i_\mathcal{V}}\cdot \Phi_{{\sigmav},ts}^{K,L}=
\int_K \dfrac{e^l-e^0}{\Delta t}\; dx + \int_{\partial K}\dfrac{p^{(l)} \bvv^{(l)} +p^{(0)}\bvv^{(0)}}{2}\cdot \bn \; d\sigma,
\end{split}
\end{equation*}
which leads to \eqref{Eq:conserv-discr}.